\documentclass[a4paper]{amsart}
\usepackage{amssymb}
\usepackage{xypic}
\usepackage{verbatim}

\def \xcirc{\objectmargin{0.1pc}\def\objectstyle{\sssize}\diagram \squarify<1pt>{}\circled\enddiagram}

\newtheorem{notations}{Notations}[section]
\newtheorem{theorem}{Theorem}[section]
\newtheorem{lemma}[theorem]{Lemma}
\newtheorem{proposition}[theorem]{Proposition}

\theoremstyle{definition}

\theoremstyle{remark}

\DeclareMathOperator{\BC}{\mathsf{BC}} \DeclareMathOperator{\BN}{\mathsf{BN}}
\DeclareMathOperator{\BP}{\mathsf{BP}} \DeclareMathOperator{\HS}{\mathsf{H}}
 \DeclareMathOperator{\HH}{\mathsf{HH}}
 \DeclareMathOperator{\HC}{\mathsf{HC}}
\DeclareMathOperator{\HP}{\mathsf{HP}} \DeclareMathOperator{\HN}{\mathsf{HN}}
\DeclareMathOperator{\ima}{\mathsf{Im}} \DeclareMathOperator{\ide}{\mathsf{id}}
\DeclareMathOperator{\Tot}{\mathsf{Tot}} 
 
 \DeclareMathOperator{\bZ}{\mathbb{Z}}

\newcommand{\ot}{\otimes}
\newcommand{\sub}{\subseteq}

\newcommand{\noi}{\noindent}
\newcommand{\wt}{\widetilde}
\newcommand{\wh}{\widehat}
\newcommand{\ov}{\overline}
\newcommand{\al}{\alpha}
\newcommand{\be}{\beta}
\newcommand{\ep}{\epsilon}
\newcommand{\de}{\delta}
\newcommand{\la}{\lambda}
\newcommand{\si}{\sigma}

\newcommand{\ba}{\mathbf a}
\newcommand{\bb}{\mathbf b}
\newcommand{\bc}{\mathbf c}

\newcommand{\bx}{\mathbf x}
\newcommand{\byy}{\mathbf y}

\newcommand{\cX}{\mathcal X}
\newcommand{\cY}{\mathcal Y}

\newcommand{\sR}{\scriptstyle R \textstyle}
\newcommand{\sdR}{\scriptstyle R \displaystyle}
\newcommand{\scw}{{}^{\!\!\!\! \scriptstyle{w}}}

\begin{document}

\title{Relative cyclic homology of square zero extensions}

\author{Jorge A. Guccione}
\address{Departamento de Matem\'atica\\ Facultad de Ciencias Exactas y Naturales, Pabell\'on
1 - Ciudad Universitaria\\ (1428) Buenos Aires, Argentina.} \curraddr{}
\email{vander@dm.uba.ar}
\thanks{Supported by PICT 12330, UBACYT 0294 and CONICET}

\author{Juan J. Guccione}
\address{Departamento de Matem\'atica\\ Facultad de Ciencias Exactas y Naturales\\
Pabell\'on 1 - Ciudad Universitaria\\ (1428) Buenos Aires, Argentina.} \curraddr{}
\email{jjgucci@dm.uba.ar}
\thanks{Supported by PICT 12330, UBACYT 0294 and CONICET}

\date{}

\dedicatory{}


\begin{abstract} Let $k$ be a characteristic zero field, $C$ a $k$-algebra and $M$ a square
zero two sided ideal of $C$. We obtain a new mixed complex, simpler than the canonical one,
giving the Hochschild and cyclic homologies of $C$ relative to $M$. This complex resembles the
canonical reduced mixed complex of an augmented algebra. We begin the study of our complex
showing that it has a harmonic decomposition like to the one considered by Cuntz and Quillen
for the normalized mixed complex of an algebra. We also give new proofs of two theorems of
Goodwillie, obtaining a light improvement of one of them.
\end{abstract}

\maketitle

\section{Introduction}
Let $k$ be a characteristic zero field, $C$ a $k$-algebra and $M$ a two-sided ideal of $C$. In
this work we deal with the  Hochschild, cyclic, periodic and negative homologies of $C$
relative to $M$, when $M^2=0$. Our main result is Theorem~\ref{th3.2}, in which we obtain a
double mixed complex $(\hat{X},\hat{b},\hat{d},\hat{B})$, given these homologies, whose
associated mixed complex $(\breve{X},\breve{b}, \breve{B})$ is simpler than the canonical
mixed complex of $C$ relative to $M$. We hope that $(\hat{X},\hat{b},\hat{d},\hat{B})$ be
useful to prove results about cyclic type homologies of an algebra relative to a nilpotent
ideal, by induction on the degree of nilpotence. Evidence in such sense is provided by
Theorem~\ref{th4.6}, in which we improve a result of Goodwillie. We also hope that
$(\hat{X},\hat{b},\hat{d},\hat{B})$ be a first step to obtain explicit computations of cyclic
homology groups of an algebra $C$ relative to a two sided square zero ideal.

\smallskip

The paper is organized in the following way:

\smallskip

In Section~2 we recall some well known definitions and results. Among them, the perturbation
lemma, which we will use again and again in the rest of the paper, and the definition of
double mixed complex, which we got from \cite{Co}.

\smallskip

Section~3 is devoted to establishing the main results in this paper. Since $M^2=0$, the
algebra $C$ is isomorphic to a square zero extension $E=A\ltimes_f M$, where $A=C/M$ and
$f\colon A \ot A\to M$ is a Hochschild normal $2$-cocycle. So, we can restrict our attention
to this type of algebras. In fact, $(\hat{X},\hat{b},\hat{d},\hat{B})$ can be thought as a
double mixed complex associated to the $3$-tuple $(A,M,f)$, and this association is functorial
in an evident sense. For $0\le 2w\le v$, let $X_v^w$ be the direct sum of all the tensor
products $X_0\ot \cdots \ot X_n$ such that $X_0=M$, $X_i=M$ for $w$ indices $i>0$ and
$X_i=\ov{A}$ for the other ones, where $n=v-w$ and $\ov{A} = A/k$. Let $b\colon X_v^w\to
X_{v-1}^w$ be the map given by the same formula as the Hochschild boundary map of an algebra,
where the meaning of the concatenation $x_ix_{i+1}$ of two consecutive factors in a simple
tensor means is the one given in item~(3) of Notation~\ref{not1.2}. Let $t\colon X_v^w\to
X_v^w$ be the map defined by
$$
t(x_0\ot\cdots\ot x_n)=(-1)^{in}x_i\ot \cdots \ot x_n\ot x_0\ot\cdots\ot x_{i-1},
$$
where $i$ denotes the last index such that $x_i\in M$ and let $N = \ide + t + t^2 + \cdots +
t^w$. The double mixed complex $(\hat{X},\hat{b},\hat{d},\hat{B})$ has objects $\hat{X}_v^w =
X_v^w\oplus X_{v-1}^w$. The boundary maps are given by
\[
\hat{b}(\bx,\byy) = \bigl(b(\bx) + (\ide-t) (\byy),-b(\byy)\bigr)\quad\text{and}\quad
\hat{d}(\bx,\byy) = \bigl(d(\bx),d'(\byy)\bigr),
\]
where $d,d'\colon X_v^w\to X_v^{w+1}$ are maps depending on $f$, and the Connes operator is
given by $\hat{B}(\bx,\byy) = (0,N(\bx))$. So, it resembles the reduced mixed complex of an
augmented algebra. Since the maps $t$ and $N$ satisfy
\begin{equation}
\ima(1-t) = \ker(N)\quad\text{and}\quad \ima(N)=\ker(1-t),\label{eq1}
\end{equation}
the cyclic homology of $E$ relative to $M$ is the homology of the quotient complex of
$(X,b,d)$ by the image of $\ide-t$. Indeed, this also follows from the fact that
$(\hat{X},\hat{b},\hat{d},\hat{B})$ satisfies the Connes property (\cite{C-Q}), which is
another consequence of the equalities~\eqref{eq1}. We finish the section giving a new proof of
the following celebrate theorem of Goodwillie: if $M$ is a nilpotent two-sided ideal of a
$k$-algebra $C$, then $\HP(C)=\HP(C/M)$.

\smallskip

\noi The aim of Section~4 is to show that $(\breve{X},\breve{b},\breve{B})$ has a harmonic
decomposition like the one studied in \cite{C-Q}. In order to carry out this task we need to
define a de Rham coboundary map and a Karoubi operator on $(\breve{X},\breve{b})$. Actually it
will be convenient for us to work with a new double mixed complex, namely
$(\ddot{X},\ddot{d},\ddot{b},\ddot{B})$, whose associated mixed complex is also
$(\breve{X},\breve{b},\breve{B})$. As in \cite{C-Q} the Karoubi operator $\ddot{\kappa}$ of
$(\ddot{X},\ddot{d},\ddot{b})$ commutes with $\ddot{b}$ and $\ddot{d}$ and satisfies a
polynomial equation $P_w(\kappa)$ on each $\ddot{X}_v^w$. Thus we have the harmonic
decomposition $\ddot{X} = P(\ddot{X})\oplus P^{\perp}(\ddot{X})$, where $P$ is the spectral
projection onto the generalized nullspace for $\ide - \ddot{\kappa}$ and $P^{\perp} = 1-P$.
The first component of this decomposition is $\ddot{B}$-acyclic and the second one is
$\ddot{d}$-acyclic and killed by $\ddot{B}$. Hence $(\ddot{X},\ddot{d},\ddot{b})$ has the
Connes property. We finish the section by giving two explicit descriptions of $P(\ddot{X})$
and obtaining a new expression for the connection map of the long exact sequence relating the
absolute Hochschild homologies of $A$ and $E$, with the Hochschild homology of $E$ relative to
$M$.

\smallskip

Although we had assume that $k$ is a characteristic zero field, many of the results in this
paper are valid under considerable weaker hypothesis. More precisely we can take a commutative
ring $k$ and a $k$ subalgebra $S$ of $A$ and consider the $S$-relative Hochschild, cyclic,
negative and periodic homologies. In this case we must replace $\ov{A}$ by $A/S$ and the
tensor products over $k$, that appear in all the complexes in this paper, by cyclic tensor
products over $S$ (See \cite{K1}, \cite{G-S}, \cite{K2} and \cite{Q})). All the results of
Section~3 are valid in this context, with the exception of Lemma~\ref{le3.3},
Theorem~\ref{th3.4} and Propositions~\ref{pr3.5} and~\ref{pr3.8}. If $k$ contains
$\mathbb{Q}$, then all the results in this paper are valid, except Theorem~\ref{th4.6}, and
this theorem is also valid if we also have that $S^e$ is semisimple. Finally, when $S^e$ is a
separable $k$-algebra, the relative and absolute homologies coincide, as was shown in the
above mentioned papers.

\smallskip

Next we introduce some notations that we will use throughout this paper.

\begin{notations} Let $k$ be a commutative ring, $V$ a $k$-module, $C$ a $k$-algebra and $M$
a $C$-bimodule.

\begin{enumerate}

\smallskip

\item We put $\ov{C} = C/k$ and given $x\in C$ we also let $x$ denote its class in $\ov{C}$.

\smallskip

\item We let $V^{\ot n}$ denote the $n$-fold power tensor of $V$.

\smallskip

\item Given $x_0\ot\cdots\ot x_n \in C\ot \ov{C}^{\ot n}$ and $0\le i<j\le n$, we write
$\bx_i^j = x_i\ot\cdots\ot x_j$.

\smallskip

\item For $n\ge 0$, we let $M^{\ot_{\!C}^n}$ denote the $n$-fold power tensor of $M$ over
$C$. As usual, we consider that $M^{\ot_{\!C}^0} = C$.

\smallskip

\item Given $x_1\ot_C\cdots\ot_C x_n \in M^{\ot_{\!C}^n}$ and $1\le i<j\le n$, we will write
$\ov{\bx}_i^j = x_i\ot_C\cdots\ot_C x_j$.
\end{enumerate}
\label{not1.1}
\end{notations}

\begin{notations} Let $E = A\ltimes_f M$ be a square zero extension.

\begin{enumerate}

\smallskip

\item We let $\pi_{\!A}\colon E\to A$ and $\pi_{\!M}\colon E\to M$ denote the maps defined by
$\pi_{\!A}(a,m)= a$ and $\pi_{\!M}(a,m)=m$, respectively.

\smallskip

\item We extend $f$ to $E\ot E$ writing $f(x,y) = 0$ if $x\in M$ or $y\in M$.

\smallskip

\item Given $x,y\in A\bigcup M$ we set
\[
xy = \begin{cases} \text{the product of $x$ and $y$ in $A$} &\text{if $x,y\in A$,}\\
\text{the left action of $x$ on $y$} &\text{if $x\in A$ and $y\in M$,}\\ \text{the right
action of $y$ on $x$} &\text{if $x\in M$ and $y\in A$,}\\ 0 &\text{if $x,y\in M$.}
\end{cases}
\]

\smallskip

\item For $0\le w\le n$, let $B^n_w\sub \ov{E}^{\ot n}$ be the $k$-submodule spanned by the
$n$-tensors $x_1\ot\cdots \ot x_n$ such that exactly $w$ of the $x_i$'s belong to $M$, while
the other ones belong to $\ov{A}$. To unify expressions we make the convention that $B_0^0 =
k$ and $B^n_w = 0$, for $w<0$ or $n<w$.

\smallskip

\item For $\bx_0^n\in M\ot B^n_w\bigcup A\ot B^n_{w+1}$ and $0\le l\le n$, we
define $\mu_l(\bx_0^n)$ by
$$
\mu_l(\bx_0^n) = \begin{cases} (-1)^l \bx_0^{l-1}\ot x_lx_{l+1}\ot \bx_{l+2}^n &\text{if $0\le
l< n$,}\\ (-1)^n x_nx_0\ot\bx_1^{n-1}&\text{If $l = n$.}\end{cases}
$$
Moreover we set
\begin{alignat*}{2}
&\qquad\mu_0^A(\bx_0^n) = \pi_{\!A}(x_0x_1)\ot\bx_2^n,\qquad &&\mu_n^A(\bx_0^n) =
(-1)^n\pi_{\!A}(x_nx_0)\ot\bx_1^{n-1},\\
&\qquad\mu_0^M(\bx_0^n) = \pi_{\!M}(x_0x_1)\ot\bx_2^n,\qquad &&\mu_n^M(\bx_0^n) =
(-1)^n\pi_{\!M}(x_nx_0)\ot\bx_1^{n-1}.
\end{alignat*}

\smallskip

\item For $\bx_0^n\in M\ot B^n_w\bigcup A\ot B^n_{w+1}$ and $0\le j\le n$, we
define $\mu_l(\bx_0^n)$ by
$$
F_j(\bx_0^n) = \begin{cases} (-1)^j \bx_0^{j-1}\ot f(x_j,x_{j+1})\ot \bx_{j+2}^n &\text{if
$0\le j< n$,}\\ - \bx_1^{n-1}\ot f(x_n,x_0) &\text{If $j = n$.}\end{cases}
$$

\smallskip

\item For an elementary tensor $\bx_0^n$, such that $x_i \in A\bigcup M$ for all $i$,
we let $i(\bx_0^n)$ denote the last index $i$ such that $x_i\in M$.

\smallskip

\item For an elementary tensor $\bx_0^n$ such that $x_i \in A\bigcup M$ for all $i$, we define
$t(\bx_0^n) = (-1)^{i(\bx_0^n)n}\bx_{i(\bx_0^n)}^n\ot \bx_0^{i(\bx_0^n)-1}$.

\end{enumerate}
\label{not1.2}
\end{notations}

\noi \textit{Acknowledgement.}\enspace We would like to thank our colleague, Professor
Guillermo Corti\-\~nas for his careful reading of a first version of our paper and for his
suggestions that have helped us to improve the presentation of this paper substantially.

\section{Preliminaries}
In this section we recall some well known definitions and results, and we fix some notations
that we will use in the rest of the paper. Let $C$ be a $k$-algebra.
\subsection{Double and triple complexes}
A double complex $\cX = (X,d^v,d^h)$ of $C$-modules, is a family $(X_{pq})_{p,q\in\bZ}$ of
$C$-modules, together with $C$-linear maps
\[
d^h\colon X_{pq}\to X_{p-1,q}\quad\text{and}\quad d^v\colon X_{pq}\to X_{p,q-1},
\]
such that $d^h\xcirc d^h=0$, $d^v\xcirc d^v=0$ and $d^v \xcirc d^h + d^h\xcirc d^v =0$. The
total complex of $(X,d^v,d^h)$ is the complex $\Tot(\cX) = (X,d)$, in which
\[
X_n = \prod_p X_{p,n-p}\quad\text{and}\quad d = d^v+d^h.
\]
A morphism of double complexes $f\colon (X,d^v,d^h)\to (Y,\de^v,\de^h)$ is a family of maps
$f\colon X_{pq}\to Y_{pq}$, such that $\de^v\xcirc f = f\xcirc d^v$ and $\de^h \xcirc f=
f\xcirc d^h$. The morphism from $\Tot(X,d^v,d^h)$ to $\Tot(Y,\de^v,\de^h)$ induced by $f$ will
be denoted $\Tot(f)$.

\smallskip

Similarly, one can give the notions of triple complex $\cX = (X,d^v,d^h,d^d)$ and of morphism
of triple complexes. For a triple complex $\cX$, there are three ways for constructing a
double complex by taking total complexes of double complexes. We call each one of these double
complexes a partial total complex of $\cX$. Finally, the total complex $\Tot(\cX)$ of $\cX$,
is the total complex of any of its partial total complexes. Of course, $\Tot(\cX)$ is
independently of the chosen way to construct it.
\subsection{Mixed complexes}
In this subsection we recall briefly the notion of mixed complex. For more details about this
concept we refer to \cite{Ka1} and \cite{B}.

\smallskip

A mixed complex $(X,b,B)$ is a graded $C$-module $(X_n)_{n\ge 0}$, endowed with morphisms
$b\colon X_n\to X_{n-1}$ and $B\colon X_n\to X_{n+1}$, such that
$$
b\xcirc b = 0,\quad B\xcirc B = 0\quad\text{and}\quad B \xcirc b + b\xcirc B = 0.
$$
A morphism of mixed complexes $f\colon (X,b,B)\to (Y,d,D)$ is a family $f\colon X_n\to Y_n$,
such that $d\xcirc f = f\xcirc b$ and $D\xcirc f= f\xcirc B$. A mixed complex $\cX = (X,b,B)$
determines a double complex
\[
\xymatrix{\\\\ \BP(\cX)=}\qquad
\xymatrix{
& \vdots \dto^-{b} &\vdots \dto^-{b}& \vdots \dto^-{b}\\
\dots & X_2 \lto_-{B}\dto^-{b} & X_1 \lto_-{B}\dto^-{b} & X_0 \lto_-{B}\\
\dots & X_1 \lto_-{B}\dto^-{b} & X_0 \lto_-{B}\\
\dots & X_0 \lto_-{B}}
\]
By deleting the positively numbered columns we obtain a subcomplex $\BN(\cX)$ of $\BP(\cX)$.
The quotient double complex $\BP(\cX)/\BN(\cX)$ is denoted by $\BC(\cX)$. The homologies
$\HC_*(\cX)$, $\HN_*(\cX)$ and $\HP_*(\cX)$, of the total complexes of $\BC(\cX)$, $\BN(\cX)$
and $\BP(\cX)$ respectively, are called the cyclic, negative and periodic homologies of $\cX$.
The homology $\HH_*(\cX)$, of $(X,b)$, is called the Hochschild homology of $\cX$. Finally, it
is clear that a morphism $f\colon \cX\to \cY$ of mixed complexes induces a morphism from the
double complex $\BP(\cX)$ to the double complex $\BP(\cY)$.

\smallskip

Following \cite{Co} by a double mixed complex we will understand a bigraded module $X$
equipped with three $k$-linear maps of degree $\pm 1$: $\partial$ that lowers the first index
and fixes the second one, $\de$ that fixes the first index and lowers the second one, and $B$
which fixes the first index and increases the second one. These maps satisfy
\[
0 = \partial^2 = \de^2 = B^2 = \de\xcirc \partial + \partial\xcirc \de = \de \xcirc B +
B\xcirc \de = \partial \xcirc B + B\xcirc \partial.
\]
The mixed complex $(X,\de + \partial,\!B)$ associated with a double mixed complex
$(X,\de,\partial,B)$ is obtained setting $(X,\de+\partial) = \Tot(X,\de,\partial)$ and $B_n =
\bigoplus_{i+j=n} B_{ij}$. By definition, the Hochschild, cyclic, periodic and negative
homologies of $(X,\de,\partial,B)$ are the Hochschild, cyclic, periodic and negative
homologies of $(X,\de+\partial,B)$, respectively.
\subsection{The relative Hochschild and cyclic homologies}
Let $C$ be a $k$-algebra and let $(C\ot \ov{C}^{\ot *},b,B)$ be the normalized mixed complex
of $C$. Recall that the cyclic, negative, periodic and Hochschild homologies $\HC_*(C)$,
$\HN_*(C)$, $\HP_*(C)$ and $\HH_*(C)$ of $C$ are the respective homologies of
$(C\ot\ov{C}^{\ot *},b,B)$.

\smallskip

Next, we define the relative homologies. Let $I$ be a two sided ideal of $C$ and let $D =
C/I$. The cyclic, negative, periodic and Hochschild homologies $\HC_*(C,I)$, $\HN_*(C,I)$,
$\HP_*(C,I)$ and $\HH_*(C,I)$, of $C$ relative to $I$, are by definition the respective
homologies of the mixed complex
\[
\xymatrix{{\ker\bigl((C\ot \ov{C}^{\ot *},b,B)} \rto^-{\pi} &{(D\ot \ov{D}^{\ot *},b,B)
\bigr)}},
\]
where $\pi$ is the map induced by the canonical projection from $C$ onto $D$.
\subsection{The perturbation lemma}
Next, we recall the perturbation lemma. We give the more general version introduced in
\cite{C}.

\smallskip

A homotopy equivalence data
\begin{equation}
\xymatrix{(Y,\partial)\ar@<-1ex>[r]_-{i} & (X,d) \ar@<-1ex>[l]_-{p}}, \quad h\colon X_*\to
X_{*+1},\label{eq2}
\end{equation}
consists of the following:

\begin{enumerate}

\smallskip

\item Chain complexes $(Y,\partial)$, $(X,d)$ and quasi-isomorphisms $i$,
$p$ between them,

\smallskip

\item A homotopy $h$ from $i\xcirc p$ to $\ide$.
\end{enumerate}

\smallskip

\noi A perturbation~$\de$ of~\eqref{eq2} is a map $\de\colon X_*\to X_{*-1}$ such that
$(d+\de)^2 = 0$. We call it small if $\ide - \de\xcirc h$ is invertible. In this case we write
$A = (\ide - \de\xcirc h)^{-1}\xcirc \de$ and we consider
\begin{equation}
\xymatrix{(Y,\partial^1)\ar@<-1ex>[r]_-{i^1} & (X,d+\de)\ar@<-1ex>[l]_-{p^1}}, \quad h^1\colon
X_*\to X_{*+1},\label{eq3}
\end{equation}
with
$$
\partial^1 = \partial + p\xcirc A\xcirc i,\quad i^1 = i + h\xcirc A\xcirc i,\quad
p^1 = p + p\xcirc A\xcirc h,\quad h^1 = h + h\xcirc A\xcirc h.
$$
A deformation retract is a homotopy equivalence data such that $p\xcirc i = \ide$. A
deformation retract is called special if $h\xcirc i = 0$, $p\xcirc h = 0$ and $h\xcirc h = 0$.

\smallskip

In all the cases considered in this paper the map $\de\xcirc h$ is locally nilpotent, and so
$(\ide - \de\xcirc h)^{-1} = \sum_{n=0}^{\infty} (\de\xcirc h)^n$.

\begin{theorem} (\cite{C}) If $\de$ is a small perturbation of the homotopy equivalence
data~\eqref{eq2}, then the perturbed data~\eqref{eq3} is a homotopy equivalence. Moreover, if
\eqref{eq2} is a special deformation retract, then~\eqref{eq3} is also.\label{th2.1}
\end{theorem}
\subsection{The suspension} The suspension of a chain complex $(X,d)$ is the
complex $(X,d)[1] = (X[1],d[1])$, defined by $X[1]_* = X_{*-1}$ and $d[1]_* = -d_{*-1}$.

\section{The relative cyclic homology of a square zero extension}
Let $A$ be a $k$-algebra, $M$ an $A$-bimodule and $f\colon A\ot A\to M$ a Hochschild normal
$2$-cocycle. The square zero extension $E = A\ltimes_f M$, of $A$ by $M$ associated with $f$,
is the direct sum $A\oplus M$ with the associative algebra structure given by
\[
(a,m)(a',m') = (aa',am' + ma' + f(a,a')).
\]
Let $C$ be a $k$-algebra and let $M$ be a two sided ideal of $C$ such that $M^2 = 0$. It is
well known that $C$ is isomorphic to a square zero extension $E$ of $A$ by $M$. In this
section we obtain a double mixed complex, simpler than the canonical one, giving the
Hochschild, cyclic, periodic and negative homologies of $E$ relative to $M$. Then we show that
the cyclic homology of $E$ relative to $M$, is also given by a still simpler complex. Finally,
we obtain convenient expressions for the connection map of the long exact sequences in
Hochschild and cyclic homologies, associated with the short exact sequence of mixed complexes
\begin{equation}
%
\xymatrix{0\rto &\ker(\pi) \rto^-{i} &{(E\ot\ov{E}^{\ot *},b,B)}\rto^-{\pi} &
{(A\ot\ov{A}^{\ot *},b,B)}\rto & 0.}\label{eq4}
\end{equation}

\subsection{Complexes for the relative Hochschild, cyclic, periodic and negative
homologies}
For $w\ge 0$ and $v\ge 2w$, let $X_v^w = M\ot B_w^{v-w}$. By convenience we put $X_v^w = 0$,
otherwise. From now on we often will use the indices $v,w$ and $n$, which always will satisfy
the relation $n=v-w$. Consider the triple diagram
\[
\def\objectstyle{\ssize}
\def\labelstyle{\ssize}
\grow{\xymatrix @! {\\\\\\\txt{$\cX = $\,\,\,}}}
\grow{\xymatrix @!R @R=-0.35pc @!C @C=0.4pc @L=0.4pt
{&&&&& \vdots\ar@{->}[3,0]_-{-b}&&&\vdots \ar@{->}[3,0]_-{b}&&&  \vdots \ar@{->}[3,0]_-{-b}\\\\
&&&&\vdots \ar@{->}[3,0]_-{-b}&&& \vdots \ar@{->}[3,0]_-{b}&&&\vdots\ar@{->}[3,0]_-{-b}\\
&&\cdots&&&{X_4^2}\ar@{->}'[0,-1] [0,-3]_(0.35){\ide-t} &&&{X_4^2}\ar@{->}'[0,-1]
[0,-3]_(0.23){N}&&&{X_4^2}
\ar@{->}'[0,-1] [0,-3]_(0.35){\ide-t}&&&\cdots \ar@{->} [0,-3]_-{N}\\
&&& \vdots \ar@{->}[3,0]_-{-b} &&& \vdots \ar@{->}[3,0]_-{b}&&& \vdots\ar@{->}[3,0]_-{-b}&\\
&\cdots &&&{X_4^1}\ar@{->}[-2,1]^-{d'} \ar@{->}'[2,0]_(0.70){-b}[3,0] \ar@{->}'[0,-1]
[0,-3]_(0.35){\ide-t} &&&{X_4^1}\ar@{->}[-2,1]^-{d} \ar@{->}'[2,0]_(0.70){b}[3,0]
\ar@{->}'[0,-1] [0,-3]_(0.23){N} &&&{X_4^1}\ar@{->}[-2,1]^-{d'}
\ar@{->}'[2,0]_(0.70){-b}[3,0]\ar@{->}'[0,-1] [0,-3]_(0.30){\ide-t}&&& \cdots \ar@{->}
[0,-3]_-{N}\\\\
\cdots &&& {X_4^0}\ar@{->}[0,-3]_(0.5){\ide-t} \ar@{->}[3,0]_-{-b} \ar@{->}[-2,1]^-{d'} &&&
{X_4^0} \ar@{->}[0,-3]_(0.4){N} \ar@{->}[3,0]_-{b} \ar@{->}[-2,1]^-{d} &&& {X_4^0}
\ar@{->}[0,-3]_(0.4){\ide-t} \ar@{->}[3,0]_-{-b} \ar@{->}[-2,1]^-{d'} &&& \cdots
\ar@{->}[0,-3]_(0.4){N}\\
&\cdots&&&{X_3^1}\ar@{->}'[2,0]_(0.70){-b}[3,0] \ar@{->}'[0,-1] [0,-3]_(0.35){\ide-t}
&&&{X_3^1} \ar@{->}'[2,0]_(0.70){b}[3,0] \ar@{->}'[0,-1] [0,-3]_(0.23){N} &&&{X_3^1}
\ar@{->}'[2,0]_(0.70){-b}[3,0]\ar@{->}'[0,-1] [0,-3]_(0.30){\ide-t}&&&
\cdots\ar@{->} [0,-3]_-{N}\\\\
\cdots &&& {X_3^0} \ar@{->}[0,-3]_(0.5){\ide-t} \ar@{->}[3,0]_-{-b} \ar@{->}[-2,1]^-{d'} &&&
{X_3^0} \ar@{->}[0,-3]_(0.4){N} \ar@{->}[3,0]_-{b} \ar@{->}[-2,1]^-{d} &&&
{X_3^0}\ar@{->}[0,-3]_(0.4){\ide-t} \ar@{->}[3,0]_-{-b}
\ar@{->}[-2,1]^-{d'}&&& \cdots \ar@{->}[0,-3]_(0.4){N}\\
&\cdots &&&{X_2^1}\ar@{->}'[0,-1][0,-3]_(0.35){\ide-t} &&&{X_2^1} \ar@{->}'[0,-1]
[0,-3]_(0.23){N} &&&{X_2^1}\ar@{->}'[0,-1] [0,-3]_(0.30){\ide-t}
&&&\cdots \ar@{->} [0,-3]_-{N}\\\\
\cdots &&& {X_2^0} \ar@{->}[0,-3]_(0.5){\ide-t} \ar@{->}[3,0]_-{-b} \ar@{->}[-2,1]^-{d'} &&&
{X_2^0} \ar@{->}[0,-3]_(0.4){N} \ar@{->}[3,0]_-{b} \ar@{->}[-2,1]^-{d} &&& {X_2^0}
\ar@{->}[0,-3]_(0.4){\ide-t} \ar@{->}[3,0]_-{-b} \ar@{->}[-2,1]^-{d'} &&& \cdots
\ar@{->}[0,-3]_(0.4){N}\\\\\\
\cdots &&& {X_1^0} \ar@{->}[0,-3]_(0.5){\ide-t} \ar@{->}[3,0]_-{-b} &&& {X_1^0}
\ar@{->}[0,-3]_(0.4){N} \ar@{->}[3,0]_-{b} &&& {X_1^0}\ar@{->}[0,-3]_(0.4){\ide-t}
\ar@{->}[3,0]_-{-b}&&& \cdots \ar@{->}[0,-3]_(0.4){N}\\\\\\
\cdots &&& {X_0^0} \ar@{->}[0,-3]_(0.5){\ide-t} &&& {X_0^0}\ar@{->}[0,-3]_(0.4){N} &&& {X_0^0}
\ar@{->}[0,-3]_(0.4){\ide-t} &&& \cdots \ar@{->}[0,-3]_(0.4){N} }}
\]
where
\begin{alignat*}{2}
& b(\bx_0^n) = \sum_{j=0}^n \mu_j(\bx_0^n),&&\qquad d(\bx_0^n) = \sum_{j=1}^{n-1}
F_j(\bx_0^n),\\
&d'(\bx_0^n) = -d(\bx_0^n) - \!\! \sum_{j=i(\bx_0^n)+1}^{n-1}\!\! t(F_j(\bx_0^n)), &&\qquad
N(\bx_0^n) = \sum_{l=0}^w t^l(\bx_0^n),
\end{alignat*}
the middle face $(X,b,d)$ is the $0$-th face and the bottom row is the $0$-th row. Note that
$\ide - t\colon X_v^0\to X_v^0$ is the zero map and $N\colon X_v^0\to X_v^0$ is the identity
map.

\smallskip

For $l\in\bZ$, let $\tau^l(\cX)$ be the subdiagram of $\cX$ obtaining by deleting the $u$-th
faces $(X,b,d)$ and $(X,-b,d')$ with $u>l$, and let $\tau_0(\cX)$ and $\tau_0^1(\cX)$ be the
quotient triple diagrams of $\cX$ by $\tau^{-1}(\cX)$ and $\tau^1(\cX)$ by $\tau^{-1}(\cX)$,
respectively.

\begin{theorem} $\cX$ is a triple complex. Moreover
\begin{align*}
&\HH_*(E,M) = \HS_*(\Tot(\tau_0^1(\cX))),\\
&\HC_*(E,M) = \HS_*(\Tot(\tau_0(\cX))),\\
&\HP_*(E,M) = \HS_*(\Tot(\cX)),\\
&\HN_*(E,M) = \HS_*(\Tot(\tau^1(\cX))).
\end{align*}
\label{th3.1}
\end{theorem}

\noi Consequently, the following equalities hold:
\begin{alignat*}{3}
& d\xcirc b = - b\xcirc d,&&\qquad d'\xcirc b = - b\xcirc d',&&\qquad d'\xcirc N = - N\xcirc
d,\\
& d\xcirc (\ide-t) = -(\ide-t)\xcirc d',&&\qquad b\xcirc N = N\xcirc b,&&\qquad t\xcirc b =
b\xcirc t.
\end{alignat*}
We will use these equalities freely throughout the paper.

\smallskip

Theorem~\ref{th3.1} is a consequence of Theorem~\ref{th3.2}, which we enounce below and whose
proof will be relegated to Appendix~A. For $w\ge 0$ and $v\ge 2w$, let $\hat{X}^w_v =
X^w_v\oplus X^w_{v-1}$. Consider the diagram $(\hat{X},\hat{b},\hat{d})$ where the maps
$\hat{b}\colon \hat{X}^w_v \to \hat{X}^w_{v-1}$ and $\hat{d}\colon \hat{X}^w_v \to
\hat{X}^{w+1}_v$ are defined by
\[
\hat{b}(\bx,\byy) = \bigl(b(\bx) + (\ide-t)(\byy),-b(\byy)\bigr)\quad\text{and}\quad
\hat{d}(\bx,\byy) = \bigl(d(\bx),d'(\byy)\bigr).
\]
Note that $(\hat{X},\hat{b},\hat{d})$ is one of the partial total complexes of the triple
complex $\tau_0^1(\cX)$, and so $\Tot(\tau_0^1(\cX)) = \Tot(\hat{X},\hat{b},\hat{d})$ (but we
have not proved that $\tau_0^1(\cX)$ is a triple complex, yet). Let $\hat{B}\colon
\hat{X}_v^w\to \hat{X}_{v+1}^w$ be the map defined by $\hat{B}(\bx,\byy) = (0,N(\bx))$.

\begin{theorem} The following assertions hold:

\begin{enumerate}

\smallskip

\item $(\hat{X},\hat{b},\hat{d},\hat{B})$ is a double mixed complex.

\smallskip

\item The Hochschild, cyclic, periodic and negative homologies of $(\hat{X},\hat{b},\hat{d},
\hat{B})$ are the Hochschild, cyclic, periodic and negative homologies of $E$ relative to $M$,
respectively.
\end{enumerate}
\label{th3.2}
\end{theorem}

\begin{proof} See Appendix~A.
\end{proof}

\noi{\bf Proof of Theorem~\ref{th3.1}}\enspace Let $(\breve{X},\breve{b},\breve{B})$ be the
mixed complex associated with $(\hat{X},\hat{b},\hat{d},\hat{B})$. Theorem~\ref{th3.1} follows
immediately from Theorem~\ref{th3.2} and the fact that
\begin{align*}
&\Tot(\tau_0^1(\cX)) = (\breve{X},\breve{b}),\\
&\Tot(\tau_0(\cX)) = \Tot(\BC(\breve{X},\breve{b},\breve{B})),\\
&\Tot(\cX) = \Tot(\BP(\breve{X},\breve{b},\breve{B})),\\
&\Tot(\tau^1(\cX)) = \Tot(\BN(\breve{X},\breve{b},\breve{B})),
\end{align*}
which can be easily checked.\qed

\begin{lemma} The rows of the triple complex $\cX$ are contractible.\label{le3.3}
\end{lemma}

\begin{proof} For $w\ge 0$ and $v\ge 2w$, let $\si,\si'\colon X^w_v\to X^w_v$ be the maps
defined by $\si = \frac{1}{w+1}\ide$ and $\si' = \sum_{j=0}^{w-1} \frac{w-j}{w+1} t^j$. A
direct computation shows that: \allowdisplaybreaks
\begin{align*}
&\si\xcirc N = N\xcirc\si = \frac{1}{w+1} N,\\
& (\ide-t)\xcirc \si' = \si'\xcirc (\ide-t) = \sum_{j=0}^{w-1} \frac{w-j}{w+1} t^j
-\sum_{j=1}^w\frac{w-j+1}{w+1} t^j = \ide - \frac{1}{w+1}N.
\end{align*}
The result follows immediately from these equalities.
\end{proof}

\begin{theorem} Let $(\ov{X},\ov{b},\ov{d})$ be the cokernel of $\ide-t\colon
(X,-b,d') \to (X,b,d)$. The relative cyclic homology $\HC_*(E,M)$ is the homology of
$(\ov{X},\ov{b},\ov{d})$.\label{th3.4}
\end{theorem}

\begin{proof} This follows immediately from Theorem~\ref{th3.1} and Lemma~\ref{le3.3}.
\end{proof}

From Lemma~\ref{le3.3} it follows also that $(\hat{X},\hat{b},\hat{d},\hat{B})$ has the Connes
property. Theorem~\ref{th3.4} can be alternatively deduced from this fact.

\smallskip

Let $C$ be a $k$-algebra, $M$ a two sided ideal of $C$ and $A=C/M$. Next we give an
alternative proof of a celebrated theorem of Goodwillie obtained in \cite{G2}. To do this we
need the following result:

\begin{proposition} It is true that $b\xcirc \si' = \si'\xcirc b$.\label{pr3.5}
\end{proposition}

\begin{proof} Fix $w\ge 0$ and $v>2w$ and let $\bx_0^n \in X_v^w$ be an elementary tensor.
Let $0=i_0<i_1<\dots<i_w\le n$ be the indexes such that $x_{i_j}\in M$ and set $i_{w+1} =
n+1$. A direct computation shows that
$$
\mu_l\xcirc t^j(\bx_0^n) = \begin{cases} t^j\xcirc \mu_{l+i_{w+1-j}}(\bx_0^n) & \text{if $0\le
l\le n-i_{w+1-j}$,}\\ t^j\xcirc \mu_{l-n-1+i_{w+1-j}}(\bx_0^n) & \text{if $n-i_{w+1-j}<l\le
n$.}
\end{cases}
$$
Hence, the families $\mu_l\xcirc t^j(\bx_0^n)$ and $t^j\xcirc \mu_l(\bx_0^n)$ coincide, and so
$$
b\xcirc \si'(\bx_0^n) = \sum_{l=0}^n\sum_{j=0}^{w-1}\frac{w-j}{w+1}\mu_l\xcirc t^j(\bx_0^n) =
\sum_{j=0}^{w-1}\sum_{l=0}^n\frac{w-j}{w+1} t^j\xcirc \mu_l(\bx_0^n) = \si'\xcirc b(\bx_0^n),
$$
as we want.
\end{proof}

\begin{theorem} (Goodwillie) Let $C$ be a $k$-algebra, $M$ a two-sided ideal of $C$ and
$A=C/M$. If $M$ is nilpotent, then $\HP_*(C) = \HP_*(A)$.\label{th3.6}
\end{theorem}

\begin{proof} Without loss of generality we can assume that $M^2 = 0$. So $C$ is isomorphic
to a squared zero extension $E = A\ltimes_f M$. We prove the theorem showing that $\HP_*(E,M)
= 0$. Let us consider the filtration $(F_q(\cX))_{q\ge 0}$ of $\cX$, given by $F_q(X^w_v) =
X^w_v$ if $w\ge q$ and $F_q(X^w_v) =  0$ otherwise. Taking the graded complex associated with
this filtration we can assume that $E$ the zero squared extension $E = A\ltimes M$, with
trivial cocycle. In this case $\cX$ is contractible, since $d=d'=0$ and, by
Proposition~\ref{pr3.5}, $b\xcirc \si' = \si'\xcirc b$.
\end{proof}

\subsection{The connection map for the Hochschild homology}\label{TheConnection}
Let
\[
\spreaddiagramcolumns{-1.18pc}
\xymatrix{\dots\rto &{\HH_n(E,M)}\rto &{\HH_n(E)}\rto &{\HH_n(A)} \rto^-{\de_n} &
{\HH_{n-1}(E,M)}\rto &{\HH_{n-1}(E)}\rto &\dots}
\]
be the long exact sequence associated with the short exact sequence~\eqref{eq4}. In this
subsection we obtain a morphism of complexes $\breve{\de}\colon (A\ot\ov{A}^{\ot *},b)\to
(\breve{X},\breve{b})[1]$, inducing the maps $\de_n$. For each $n\ge 0$, let
\[
\de^1\colon A\ot \ov{A}^{\ot n}\to X_{n-1}^0,\quad \de^2\colon A\ot \ov{A}^{\ot n}\to
X_{n-2}^0\quad \text{and} \quad \de^3\colon A\ot \ov{A}^{\ot n}\to X_{n-1}^1
\]
be the maps defined by
$$
\de^1= \sum_{j=0}^n t\xcirc F_j,\quad \de^2 = \mu_0\xcirc F_1\quad\text{and}\quad \de^3 =
\sum_{i=2}^{n-1} F_0\xcirc F_i - \sum_{0\le i<j\le n} t\xcirc F_i\xcirc F_j.
$$
The proof of Proposition~\ref{pr3.7} below and Proposition~\ref{pr3.8} in the next subsection
are relegated to Appendix~A.

\begin{proposition} The connection map $\de_n\colon \HH_n(A)\to \HH_{n-1}(E,M)$ is induced by
the morphism of complexes $\breve{\de}\colon (A\ot\ov{A}^{\ot *},b)\to
(\breve{X},\breve{b})[1]$, given by
\[
\breve{\de}(\ba) = (\de^2(\ba),\de^1(\ba),\de^3(\ba),0,0,\dots),
\]
where we are writing
\[
\breve{X}_{n-1} = X^0_{n-2}\oplus X^0_{n-1}\oplus X^1_{n-1}\oplus X^1_n \oplus X^2_n\oplus
X^2_{n+1}\oplus X^3_{n+1}\oplus X^3_{n+2}\oplus\cdots
\]
\label{pr3.7}
\end{proposition}

\begin{proof} See Appendix~A.
\end{proof}

\subsection{The connection map for the cyclic homology}
Let
\[
\spreaddiagramcolumns{-1.18pc}
\xymatrix{\dots\rto &{\HC_n(E,M)}\rto &{\HC_n(E)}\rto &{\HC_n(A)} \rto^-{\de_n} &
{\HC_{n-1}(E,M)} \rto &{\HC_{n-1}(E)}\rto &\dots}
\]
be the long exact sequence associated with the short exact sequence~\eqref{eq4}.

\begin{proposition} The connection map $\de_n\colon \HC_n(A) \to \HC_{n-1}(E,M)$ is induced
by the morphism of complexes $\ov{\de}\colon \Tot(\BC(A\ot\ov{A}^{\ot *},b,B)) \to
\Tot(\ov{X},\ov{b},\ov{d})[1]$, given by $\ov{\de}(\ba_0^n,\bb_0^{n-2},\bc_0^{n-4},\dots) =
\de^1(\ba_0^n)$. Note that the image of $\ov{\de}$ is included in
$\ov{X}^0_{n-1}$.\label{pr3.8}
\end{proposition}

\begin{proof} See Appendix~A.
\end{proof}

\section{The harmonic decomposition}
As in the proof of Theorem~\ref{th3.1}, we let $(\breve{X},\breve{b},\breve{B})$ denote the
mixed complex associated with the double mixed complex $(\hat{X},\hat{b},\hat{d},\hat{B})$,
introduced in Theorem~\ref{th3.2}. The aim of this section is to show that
$(\breve{X},\breve{b},\breve{B})$ has a harmonic decomposition like the one studied in
\cite{C-Q}. In order to carry out this task we need to define a de Rham coboundary map and a
Karoubi operator on $(\breve{X},\breve{b})$. As we said in the introduction we are going to
work with a new double mixed complex $(\ddot{X},\ddot{d},\ddot{b},\ddot{B})$, whose associated
mixed complex is also $(\breve{X},\breve{b},\breve{B})$. In the first three subsections we
follow closely the exposition of \cite{C-Q}.

\subsection{The Rham coboundary map and the Karoubi operator}
It is easy to see that $\tau_0^1(\cX)$ is the total complex of the double complex
\[
\spreaddiagramcolumns{0.75pc} \spreaddiagramrows{-0.25pc}
\xymatrix{\\\\\\\\
(\ddot{X},\ddot{b},\ddot{d})=}\qquad\quad
\xymatrix{\vdots \dto^-{\ddot{b}}&\vdots\dto^-{\ddot{b}}& \vdots \dto^-{\ddot{b}} & \vdots
\dto^-{\ddot{b}}\\
\ddot{X}^2_4 & \ddot{X}^1_4 \dto^{\ddot{b}}\lto_{\ddot{d}}&
\ddot{X}^0_4\dto^{\ddot{b}} \lto_{\ddot{d}} & \ddot{X}^{-1}_4\dto^{\ddot{b}}\lto_{\ddot{d}}\\
& \ddot{X}^1_3 \dto^{\ddot{b}} & \ddot{X}^0_3\dto^{\ddot{b}}\lto_{\ddot{d}} &
\ddot{X}^{-1}_3\dto^{\ddot{b}}\lto_{\ddot{d}}\\
& \ddot{X}^1_2 & \ddot{X}^0_2 \dto^{\ddot{b}}\lto_{\ddot{d}} & \ddot{X}^{-1}_2
\dto^{\ddot{b}}\lto_{\ddot{d}}\\
&& \ddot{X}^0_1 \dto^{\ddot{b}} & \ddot{X}^{-1}_1 \dto^{\ddot{b}}\lto_{\ddot{d}}\\
&& \ddot{X}^0_0 & \ddot{X}^{-1}_0 \lto_{\ddot{d}} }
\]
where $\ddot{X}^w_v = X^w_v\oplus X^{w+1}_v$ and the boundary maps are defined by
\[
\ddot{b}(\bx,\byy) = (b(\bx),-b(\byy))\quad\text{and}\quad \ddot{d}(\bx,\byy) = (d(\bx)+
(\ide-t)(\byy),d'(\byy)).
\]
The de Rham coboundary map $\ddot{d}\sR\colon \ddot{X}_v^w\to \ddot{X}_v^{w-1}$ is defined by
$\ddot{d}\sR(\bx,\byy)=(0,\bx)$. It is immediate that $(\ddot{X},\ddot{d}\sR)$ is acyclic. We
now define the Karoubi operator of $\ddot{X}$. Let $\ddot{\kappa}(0)\colon \ddot{X}_v^w\to
\ddot{X}_v^w$ and $\ddot{\kappa}(1)\colon \ddot{X}_v^w\to\ddot{X}_v^w$ be the maps defined by
\[
\ddot{\kappa}(0)(\bx,\byy) = (t(\bx),t(\byy))\quad\text{and}\quad \ddot{\kappa}(1)(\bx,\byy) =
\bigl(0,d(\bx) - d'(\bx)\bigr).
\]
The Karoubi operator $\ddot{\kappa}$ of $\ddot{X}$ is the degree zero operator defined by
\[
\ddot{\kappa} = \ddot{\kappa}(0) + \ddot{\kappa}(1).
\]
Let $\breve{d}\sR\colon \breve{X}_n\to \breve{X}_{n+1}$ and $\breve{\kappa}\colon
\breve{X}_n\to \breve{X}_n$ be the maps defined by
$$
\breve{d}\sdR_n = \bigoplus_{w=-1}^n \ddot{d}\sdR_{n+w}^{\scw}\quad\text{and}\quad
\breve{\kappa}_n = \bigoplus_{w=-1}^n \ddot{\kappa}_{n+w}^w,
$$
respectively. A direct computation shows that
\begin{equation}
\ide - \ddot{\kappa} = \ddot{d} \xcirc \ddot{d}\sdR+ \ddot{d}\sdR \xcirc \ddot{d}
\quad\text{and}\quad 0 = \ddot{b} \xcirc \ddot{d}\sdR+ \ddot{d}\sdR \xcirc \ddot{b}.
\label{eq5}
\end{equation}
In particular, $\ddot{\kappa}$ is homotopic to the identity with respect to either of the
differentials $\ddot{d}$, $\ddot{d}\sR$, and so it commutes with them. From~\eqref{eq5} it
follows that
\[
\ide - \breve{\kappa} = \breve{b} \xcirc  \breve{d}\sdR +  \breve{d}\sdR \xcirc \breve{b}.
\]
Consequently, $\breve{\kappa}$ commutes with $\breve{b}$ and $\breve{d}\sR$. Hence,
$\ddot{\kappa}$ also commutes with $\ddot{b}$ (which can be also proved by a direct
computation). Let $\ddot{B}\colon \ddot{X}_v^w\to \ddot{X}_v^{w-1}$ be the map defined by
$\ddot{B}(\bx,\byy)= (0,N(\bx))$. An easy computation shows that
$(\ddot{X},\ddot{d},\ddot{b},\ddot{B})$ is a double mixed complex and that its associated
mixed complex is $(\breve{X},\breve{b},\breve{B})$. Furthermore, $\ddot{B}(\bx) = \sum_{i=0}^w
\ddot{\kappa}^i \xcirc \ddot{d}\sR(\bx)$ for all $\bx\in \ddot{X}_v^w$. Using this we obtain:
\[
\ddot{B}\xcirc  \ddot{\kappa} =  \ddot{\kappa}\xcirc \ddot{B} = \ddot{B} \quad\text{and}\quad
\ddot{d}\sdR\xcirc \ddot{B} = \ddot{B} \xcirc  \ddot{d}\sdR = 0.
\]

\subsection{The harmonic decomposition}
From the definition of $\ddot{\kappa}$ it follows immediately that
\[
(\ddot{\kappa}^{w+2}-\ide) \xcirc (\ddot{\kappa}^{w+1}-\ide)(\ddot{X}_v^w) \subseteq
(\kappa^{w+2}-\ide)(X_{1v}^{w+1}) = 0.
\]
This implies that $\ddot{\kappa}$ satisfies the polynomial equation $P_w(\ddot{\kappa}) = 0$
on $\ddot{X}_v^w$, where
\[
P_w = (X^{w+1}-1)(X^{w+2}-1).
\]
The roots of $P_w$ are the $r$-th roots of unity, with $r = w+1$ and $r=w+2$. Moreover, $1$ is
a double root and the all other roots are simple. Consequently $\ddot{X}_v^w$ decomposes into
the direct sum of the generalized eigenspace $\ker(\ddot{\kappa} - \ide)^2$ and its complement
$\ima(\ddot{\kappa} - \ide)^2$. Combining this for all $v,w$ we obtain the following
decomposition
\[
\ddot{X} = \ker(\ddot{\kappa} - \ide)^2 \oplus\ima(\ddot{\kappa} - \ide)^2,
\]
Each of these generalized subspaces is stable under any operator commuting with
$\ddot{\kappa}$, for instance, $\ddot{b}$, $\ddot{d}$, $\ddot{d}\sR$ and $\ddot{B}$.

\subsection{The harmonic projection and the Green operator}
Let $P$ be the harmonic projection operator, which is the identity map on $\ker(\ddot{\kappa}
- \ide)^2$ and the zero map on $\ima(\ddot{\kappa} - \ide)^2$. Thus we have
\[
\ddot{X} = P(\ddot{X})\oplus P^{\perp}(\ddot{X}),
\]
where $P^{\perp} = \ide-P$. It is immediate that $(P(\ddot{X}),\ddot{d},\ddot{b},\ddot{B})$
and $(P^{\perp}(\ddot{X}),\ddot{d},\ddot{b},\ddot{B})$ are double mixed subcomplexes of
$(\ddot{X},\ddot{d},\ddot{b},\ddot{B})$. On $P^{\perp}(\ddot{X})$ the operator
$$
\ide - \ddot{\kappa} = \ddot{d} \xcirc \ddot{d}\sR + \ddot{d}\sR \xcirc \ddot{d}
$$
is both invertible and homotopic to zero with respect to either differential $\ddot{d}$ and
$\ddot{d}\sR$. Hence the complexes $(P^{\perp}(\ddot{X}),\ddot{d})$ and
$(P^{\perp}(\ddot{X}),\ddot{d}\sR)$ are acyclic. Let
\[
P(\breve{X}_n) = \bigoplus_{w=-1}^n P(\ddot{X}_{n+w}^w)\qquad\text{and}\qquad
P^{\perp}(\breve{X}_n) = \bigoplus_{w=-1}^n P^{\perp}(\ddot{X}_{n+w}^w).
\]
The same argument shows that $(P^{\perp}(\breve{X}),\breve{b})$ and
$(P^{\perp}(\breve{X}),\breve{d}\sdR)$ are also acyclic. The Green operator $G\colon
\ddot{X}\to \ddot{X}$ is defined to be zero on $P(\ddot{X})$ and the inverse of
$\ide-\ddot{\kappa}$ on $P^{\perp}(\ddot{X})$. It is immediate that
\begin{equation}
G \xcirc P = P\xcirc G = 0\quad\text{and}\quad P^{\perp} = G\xcirc (\ide-\ddot{\kappa}) =
G\xcirc (\ddot{d} \xcirc \ddot{d}\sdR + \ddot{d}\sdR \xcirc \ddot{d}).\label{eq6}
\end{equation}
Moreover $P$ and $G$ commute with each operator that commutes with $\ddot{\kappa}$.

\begin{proposition} One has
\[
P^{\perp}(\ddot{X}) = \ddot{d}\xcirc P^{\perp}(\ddot{X}) \oplus \ddot{d}\sR\xcirc
P^{\perp}(\ddot{X}).
\]
Furthermore $\ddot{d}\sR$ maps $\ddot{d}\xcirc P^{\perp}(\ddot{X})$ isomorphically onto
$\ddot{d}\sR\xcirc P^{\perp}(\ddot{X})$ with inverse $G\xcirc \ddot{d}$ and $\ddot{d}$ maps
$\ddot{d}\sR\xcirc P^{\perp}(\ddot{X})$ isomorphically onto $\ddot{d}\xcirc
P^{\perp}(\ddot{X})$ with inverse $G\xcirc \ddot{d}\sR$. This gives a new proof that
$(P^{\perp}(\ddot{X}),\ddot{d})$ and $(P^{\perp}(\ddot{X}),\ddot{d}\sR)$ are
acyclic.\label{pr4.1}
\end{proposition}

\begin{proof} The proof of Proposition~2.1 of \cite{C-Q} works in our setting.
\end{proof}

\begin{proposition} One has $\bx\in P(\ddot{X}_v^w)$ if and only if $\ddot{d}\sR(\bx)$
and $\ddot{d}\sR \xcirc \ddot{d}(\bx)$ are $\ddot{\kappa}$-invariant.\label{pr4.2}
\end{proposition}

\begin{proof} The proof of Proposition~2.2 of \cite{C-Q} works in our setting.
\end{proof}

Fix $w\ge 0$ and $v\ge 2w$. Let $\grave{X}_v^w$ and $\acute{X}_v^w$ be the image of the
canonical inclusions of $X_v^w$ into $\ddot{X}_v^w$ and $\ddot{X}_v^{w-1}$ respectively, and
let $\acute{\kappa}\colon \acute{X}_v^w\to \acute{X}_v^w$ be the map induced by
$\ddot{\kappa}$. For $\bx = (\bx_0,0)\in \grave{X}_v^w$ we write
\[
\grave{d}(\bx) = (d(\bx_0),0),
\]
and, for $\byy = (0,\byy_0)\in \acute{X}_v^w$, we write
\[
\acute{d}(\byy) = (0,d(\byy_0)),\quad \acute{d}'(\byy) = (0,d'(\byy_0))\quad\text{and}\quad
\acute{t}(\byy) = (0,t(\byy_0)).
\]
It is immediate that $\acute{\kappa}$ coincides with $\acute{t}$. Note that $\ddot{\kappa}$
has finite order on $\ddot{d}\sR(\ddot{X}) = \acute{X}$ in each degree. In fact
$\ddot{\kappa}^{w+1} = \ide$ on $\acute{X}_v^w$. By the discussion in the page~86 of
\cite{C-Q},
\begin{align}
& \acute{X}_v^w = \ker(\ide-\acute{\kappa})\oplus \ima(\ide-\acute{\kappa}),\label{eq7}\\
& P(\ddot{d}\sdR(\grave{X}_v^w)) = P(\acute{X}_v^w)=\ker(\ide-\acute{\kappa}),\\
& P^{\perp}(\ddot{d}\sdR(\grave{X}_v^w))=P^{\perp}(\acute{X}_v^w)= \ima(\ide-\acute{\kappa}),
\end{align}
and the maps $P_{\acute{X}_v^w}$ and $G_{\acute{X}_v^w}$, defined as the projection onto
$\ker(\ide-\acute{\kappa})$ associated with~\eqref{eq7} and the Green operator for
$\ide-\acute{\kappa}\colon \acute{X}_v^w\to \acute{X}_v^w$, respectively, satisfy:
\begin{equation}
P_{\acute{X}_v^w} = \frac{1}{w+1} \sum_{i=0}^w \acute{\kappa}^i\quad\text{and}\quad
G_{\acute{X}_v^w} = \frac{1}{w+1} \sum_{i=0}^w (\frac{w}{2}-i)\acute{\kappa}^i.\label{eq10}
\end{equation}
Consequently, for all $\bx\in \ddot{X}_v^w$,
\begin{align}
& P\xcirc \ddot{d}\sdR(\bx) = \frac{1}{w+1}\sum_{i=0}^w \ddot{\kappa}^i\xcirc
\ddot{d}\sdR(\bx) = \frac{1}{w+1} \ddot{B}(\bx),\label{eq11}\\
& G\xcirc \ddot{d}\sdR(\bx) = \frac{1}{w+1}\sum_{i=0}^w (\frac{w}{2}-i) \ddot{\kappa}^i\xcirc
\ddot{d}\sdR(\bx).\label{eq12}
\end{align}
The formula~\eqref{eq11} has the following consequences: it implies that
\begin{equation}
\ddot{B}(P^{\perp}(\ddot{X}))=0.\label{eq13}
\end{equation}
Using this, we obtain that $\ddot{B}(\bx) = (w+1) \ddot{d}\sdR(P(\bx))$ for all $\bx\in
\ddot{X}_v^w$. So, since $(P^{\perp}(\ddot{X}),\ddot{d}\sR)$ is acyclic,
\begin{equation}
H_*(P(\ddot{X}),\ddot{B}) = H_*(P(\ddot{X}),\ddot{d}\sdR) = H_*(\ddot{X},\ddot{d}\sdR) = 0.
\label{eq14}
\end{equation}
In the terminology of \cite{C-Q} this says that $(P(\ddot{X}),\ddot{d},\ddot{b}, \ddot{B})$ is
$\ddot{B}$-acyclic. Lastly,~\eqref{eq10} combined with~\eqref{eq12} and the second formula
of~\eqref{eq6}, allows us to obtain an explicit formula for $P$. In fact, for $\bx\in
\acute{X}_v^w$, this is given by~\eqref{eq10}. Then, assume that $\bx\in\grave{X}_v^w$. Since
by~\eqref{eq12}, $G \xcirc \ddot{d}\sdR(\bx)\in \acute{X}_v^w$, we have:
\[
G\xcirc \ddot{d}\xcirc \ddot{d}\sdR(\bx) = \ddot{d}\xcirc G\xcirc \ddot{d}\sdR(\bx) =
\acute{d}' \xcirc G \xcirc \ddot{d}\sdR(\bx) + sw\xcirc (\ide - \acute{t})\xcirc G \xcirc
\ddot{d}\sdR(\bx),
\]
where $sw\colon \acute{X}_v^w\to \grave{X}_v^w$ is the map defined by $sw(0,\bx) = (\bx,0)$.
Using this, the second formula of~\eqref{eq6}, and the fact that $\acute{t}\xcirc
\ddot{\kappa}^i \xcirc \ddot{d}\sR(\bx) = \ddot{\kappa}^{i+1}\xcirc \ddot{d}\sR(\bx)$, we
obtain:
\begin{align*}
P(\bx) & =  \bx -  G\xcirc \ddot{d}\sdR\xcirc \grave{d}_0(\bx) - G\xcirc \ddot{d}\xcirc
\ddot{d}\sdR(\bx)\\
& = \bx - \frac{1}{w+2}\sum_{i=0}^{w+1} (\frac{w+1}{2}-i) \ddot{\kappa}^i\xcirc \ddot{d}\sdR
\xcirc \grave{d}(\bx) + \frac{1}{w+1}\sum_{i=0}^w (\frac{w}{2}-i) \acute{d}'\xcirc
\ddot{\kappa}^i \xcirc \ddot{d}\sdR(\bx)\\
& - \frac{1}{w+1}\sum_{i=0}^w (\frac{w}{2}-i)sw\xcirc(\ide-\acute{t})\xcirc \ddot{\kappa}^i
\xcirc \ddot{d}\sdR(\bx)\\
& = \frac{1}{w+1}sw\xcirc\ddot{B}_n(\bx) - \frac{1}{w+2}\sum_{i=0}^{w+1}
(\frac{w+1}{2}-i) \ddot{\kappa}^i\xcirc \ddot{d}\sdR \xcirc \grave{d}(\bx)\\
&+ \frac{1}{w+1}\sum_{i=0}^w (\frac{w}{2}-i) \acute{d}'\xcirc \ddot{\kappa}^i \xcirc
\ddot{d}\sdR(\bx).
\end{align*}
We now consider the chain complex $(\ddot{X},\ddot{b},\ddot{d})$ and denote by
$\ker(\ddot{B})$, $\ima(\ddot{B})$ the kernel and image of $\ddot{B}$ on $\ddot{X}$. These are
subcomplexes of $(\ddot{X},\ddot{b},\ddot{d})$. By~\eqref{eq13} and~\eqref{eq14}, we have
$\ker(\ddot{B})/\ima(\ddot{B}) = P^{\perp}(\ddot{X})$. Consequently,
$$
H_*(\ker(\ddot{B})/\ima(\ddot{B}),\ddot{b},\ddot{d}) = 0.
$$
That is, the double mixed complex $(\ddot{X},\ddot{d},\ddot{b},\ddot{B})$ has the Connes
property (\cite{C-Q}).

\smallskip

Let us define the reduced cyclic complex $\ov{C}_X^{\la}$ to be the quotient double complex
$\ov{C}_X^{\la} = \ddot{X}/\ker(\ddot{B})$. It is easy to check that $\ov{C}_X^{\la} =
\frac{P(\ddot{X})\oplus P^{\perp}(\ddot{X})}{\ima(\ddot{B})\oplus P^{\perp}(\ddot{X})} =
\frac{P(\ddot{X})}{\ima(\ddot{B})}$ and that $\breve{B}$ induces the isomorphism of complexes
$\Tot(\ov{C}_X^{\la})[1] \simeq \ima(\breve{B})$. So, we have a short exact sequence of double
complexes
\[
\xymatrix{0 \rto &{\Tot(\ov{C}_X^{\la})[1]}\rto^-{i} & {\Tot(P)(\breve{X})}\rto^-{j}
&{\Tot(\ov{C}_X^{\la})} \rto & 0,}
\]
where $j$ is the canonical surjection and $i$ is induced by $\ddot{B}$.

\subsection{A description of $P(\ddot{X})$}\label{sAdescription}
The aim of this subsection is to obtain a precise description of the double mixed complex
$(P(\ddot{X}),\ddot{d},\ddot{b},\ddot{B})$. The main result are Theorem~\ref{th4.3} and
Proposition~\ref{pr4.4}. We relegate their proofs to Appendix~B.

\smallskip

Take $\bx = (\bx_0,\bx_1)\in \ddot{X}_v^w$, with $\bx_0\in X_v^w$ and $\bx_1\in X_v^{w+1}$. By
Proposition~\ref{pr4.2} we know that $\bx\in P(\ddot{X})$ if and only if $\bx_0$ and
$d(\bx_0)+(\ide-t)(\bx_1)$ are $t$-invariant. From this it follows immediately that if $\bx\in
P(\ddot{X})$, then $(\bx_0,\bx'_1)\in P(\ddot{X})$ for all $\bx'_1\in X_v^{w+1}$ such that
$\bx'_1-\bx_1$ is $t$-invariant. Conversely, if $\bx$ and $(\bx_0,\bx'_1)$ belong to
$P(\ddot{X})$, then $(\ide-t)(\bx'_1 - \bx_1)$ is $t$-invariant, but this implies that
$(\ide-t)(\bx'_1-\bx_1) = 0$. In other words, that $\bx'_1-\bx_1$ is $t$-invariant. In
particular $(0,\bx_1)\in P(\ddot{X})$ if and only if $\bx_1$ is $t$-invariant. We let ${}^t\!
P(\ddot{X}_v^w)$ denote the set of all elements of the shape $(0,\bx_1)\in \ddot{X}_v^w$ with
$\bx_1$ a $t$-invariant element. It is immediate that $({}^t\!P(\ddot{X}),-b,d')$ is a
subcomplex of $(P(\ddot{X}),\ddot{b},\ddot{d})$.

We assert that if $\bx\in X_v^w$ is $t$-invariant, then
\begin{equation}
(\bx,-\si'\xcirc d(\bx))\in P(\ddot{X}).\label{eq15}
\end{equation}
This follows immediately from the equality $(\ide-t)\xcirc \si'(\byy) = \byy - \frac{1}{w+2}
N(\byy)$, which was obtained in the proof of Lemma~\ref{le3.3}.

\smallskip

Recall from Theorem~\ref{th3.4} that $\ov{X}_v^w$ is the cokernel of $\ide-t\colon X_v^w\to
X_v^w$. Let $\wt{X}^w_v = \ov{X}^w_v\oplus \ov{X}^{w+1}_v$. Consider the diagram
\[
\spreaddiagramcolumns{0.75pc} \spreaddiagramrows{-0.25pc}
\xymatrix{\\\\\\\\
(\wt{X},\wt{\mathfrak{b}},\wt{\mathfrak{d}})=}\qquad\quad
\xymatrix{\vdots \dto^-{\wt{\mathfrak{b}}}&\vdots\dto^-{\wt{\mathfrak{b}}}& \vdots
\dto^-{\wt{\mathfrak{b}}} & \vdots \dto^-{\wt{\mathfrak{b}}}\\
\wt{X}^2_4 & \wt{X}^1_4 \dto^{\wt{\mathfrak{b}}}\lto_{\wt{\mathfrak{d}}}& \wt{X}^0_4
\dto^{\wt{\mathfrak{b}}} \lto_{\wt{\mathfrak{d}}} & \wt{X}^{-1}_4\dto^{\wt{\mathfrak{b}}}
\lto_{\wt{\mathfrak{d}}}\\
& \wt{X}^1_3 \dto^{\wt{\mathfrak{b}}} & \wt{X}^0_3\dto^{\wt{\mathfrak{b}}}
\lto_{\wt{\mathfrak{d}}} & \wt{X}^{-1}_3\dto^{\wt{\mathfrak{b}}}\lto_{\wt{\mathfrak{d}}}\\
& \wt{X}^1_2 & \wt{X}^0_2 \dto^{\wt{\mathfrak{b}}}\lto_{\wt{\mathfrak{d}}} & \wt{X}^{-1}_2
\dto^{\wt{\mathfrak{b}}}\lto_{\wt{\mathfrak{d}}}\\
&& \wt{X}^0_1 \dto^{\wt{\mathfrak{b}}} & \wt{X}^{-1}_1 \dto^{\wt{\mathfrak{b}}}
\lto_{\wt{\mathfrak{d}}}\\
&& \wt{X}^0_0 & \wt{X}^{-1}_0 \lto_{\wt{\mathfrak{d}}}}
\]
where the maps $\wt{\mathfrak{b}}\colon \wt{X}^w_v\to\wt{X}^w_{v-1}$ and $\wt{\mathfrak{d}}
\colon \wt{X}^w_v \to \wt{X}^{w+1}_v$ are defined by
\[
\qquad \quad \wt{\mathfrak{b}}(\bx,\byy) = (\ov{b}(\bx),-\ov{b}(\byy))\quad \text{and}\quad
\wt{\mathfrak{d}}(\bx,\byy) = (\ov{d}(\bx),-\ov{d}(\byy)),
\]
respectively. Let $\mathfrak{p}\colon X_v^w\to \ov{X}_v^w$ be the map defined by
\[
\mathfrak{p}(\bx) = \frac{n-i_w+2}{v+2}[\bx]\qquad\text{for each elementary tensor $\bx\in
X_v^w$,}
\]
where $i_w$ is the last index such that $x_{i_w}\in M$ and $[\bx]$ denotes the class of $\bx$
in $\ov{X}_v^w$. Let $\ov{N}\colon \ov{X}_v^w\to X_v^w$ be the map induced by $N$. It is
immediate that $\mathfrak{p}$ is a retraction of $\ov{N}$. Given a $t$-invariant element
$\bx\in X_v^w$, let
$$
\Upsilon(\bx) = (-\ide + \ov{N}\xcirc \mathfrak{p})\xcirc \si'\xcirc d(\bx)\in X_v^{w+1}.
$$
Let ${}^e\wt{\xi}\colon \wt{X}_v^w\to \wt{X}_{v-1}^w$ and ${}^e\wt{\varsigma}\colon
\wt{X}_v^w\to \wt{X}_v^{w+1}$ be the maps defined by
$$
{}^e \wt{\xi}(\bx,\byy) = (0,{}^e\ov{\xi}(\bx))\qquad\text{and}\qquad
{}^e\wt{\varsigma}(\bx,\byy) = (0,{}^e\ov{\varsigma}(\bx)),
$$
where ${}^e\ov{\xi}\colon \ov{X}_v^w\to \ov{X}_{v-1}^{w+1}$ and ${}^e\ov{\varsigma} \colon
\ov{X}_v^w\to \ov{X}_v^{w+2}$ are the maps given by
\begin{align*}
& {}^e\ov{\xi}(\bx) = -\frac{1}{w+1}\mathfrak{p}\xcirc \si'\xcirc d\xcirc \ov{N} \xcirc
\ov{b}(\bx) -\frac{1}{w+1}\ov{b}\xcirc \mathfrak{p}\xcirc \si'\xcirc d\xcirc\ov{N}(\bx)\\
\intertext{and}
&{}^e\ov{\varsigma}(\bx) = \frac{1}{w+1}\ov{d}\xcirc \mathfrak{p}\xcirc \si'\xcirc d\xcirc
\ov{N} (\bx) + \frac{1}{w+2} \mathfrak{p}\xcirc\si'\xcirc d\xcirc\ov{N}\xcirc\ov{d} (\bx),
\end{align*}
respectively.

\begin{theorem}\label{th4.3} Let $\wt{\mathfrak{B}}\colon \wt{X}_v^w\to
\wt{X}_{v+1}^w$ be the map defined by $\wt{\mathfrak{B}}(\bx,\byy) = (0,\bx)$. The diagrams
$(\wt{X},\wt{\mathfrak{d}},\wt{\mathfrak{b}} +{}^e\wt{\xi},\wt{\mathfrak{B}})$ and
$(\wt{X},\wt{\mathfrak{d}}+{}^e\wt{\varsigma},\wt{\mathfrak{b}},\wt{\mathfrak{B}})$ are mixed
double complexes and the maps
\[
\Psi\colon (\wt{X},\wt{\mathfrak{d}},\wt{\mathfrak{b}}+{}^e\wt{\xi},\wt{\mathfrak{B}})\to
(P(\ddot{X}),\ddot{d},\ddot{b},\ddot{B})\quad\text{and}\quad \Lambda\colon
(\wt{X},\wt{\mathfrak{d}},\wt{\mathfrak{b}}+{}^e\wt{\xi},\wt{\mathfrak{B}}) \to
(\wt{X},\wt{\mathfrak{d}}+{}^e\wt{\varsigma},\wt{\mathfrak{b}},\wt{\mathfrak{B}}),
\]
defined by $\Psi_v^{-1}(0,\byy) = \Lambda_v^{-1}(0,\byy) = (0,\byy)$, and
\begin{align*}
&\Psi_v^w(\bx,\byy) = \frac{1}{w+1}\left(\ov{N}(\bx),\Upsilon\xcirc \ov{N}(\bx)\right)
+ (0,\ov{N}(\byy)),\\
&\Lambda_v^w(\bx,\byy) = (\bx,\byy)+\frac{1}{w+1}(0,\mathfrak{p}\xcirc\si'\xcirc d\xcirc
\ov{N}(\bx)),
\end{align*}
for $w\ge 0$, are isomorphisms of double mixed complexes.
\end{theorem}

\begin{proof} See Appendix~B.
\end{proof}

We now give a formula for ${}^e\ov{\varsigma}$.

\begin{proposition} Let $\bx_0^n\in X_v^w$ be an elementary tensor, let $0=i_0<\dots
<i_w \le n$ be the indices such that $x_{i_j}\in M$ and let $i_{w+1} = n+1$. Given $0\le \al
\le n$ we let $j(\al)$ denote the number defined by $i_{j(\al)}\le \al < i_{j(\al)+1}$. We
have:
\[
{}^e\ov{\varsigma}([\bx_0^n]) = \sum_{\al<\be} \la_{\al\be}^{(w)} [F_{\al}\xcirc
F_{\be}(\bx_0^n)],
\]
where $[F_{\al}\xcirc F_{\be}(\bx_0^n)]$ denotes the class of $F_{\al}\xcirc F_{\be}(\bx_0^n)$
in $X_v^{w+2}$ and
\[
\la_{\al\be}^{(w)}= \frac{2(j(\be)-j(\al))} {(w+1)(w+2)(w+3)}- \frac{1}{(w+2)(w+3)}.
\]
\label{pr4.4}
\end{proposition}

\begin{proof} See Appendix~B.
\end{proof}

Let $i\colon \ov{X}_v^{w+1}\to \wt{X}_v^w$ and $\pi\colon \wt{X}_v^w\to \ov{X}_v^w$ be the
canonical maps. The short exact sequence of double complexes
\begin{equation}
\spreaddiagramcolumns{-0.5pc}
\xymatrix{0\rto &{(\ov{X}_*^{*+1},\ov{b},\ov{d})} \rto^-i
&{(\wt{X}_*^*,\wt{\mathfrak{b}},\wt{\mathfrak{d}}+{}^e\wt{\varsigma})} \rto^-{\pi} &
{(\ov{X}_*^*,\ov{b},\ov{d})}\rto & 0.}\label{eq16}
\end{equation}
splits in each level via the maps $s\colon \ov{X}_v^w\to\wt{X}_v^w$ and $r\colon \wt{X}_v^w
\to \ov{X}_v^{w+1}$, given by $s(\bx) = (\bx,0)$ and $r(\bx,\byy) = \byy$. From this it
follows immediately that the connection map of the homology long exact sequence associated
with~\eqref{eq16} is induced by the morphism of double complexes ${}^e\ov{\varsigma}\colon
(\ov{X}_*^*,\ov{b},\ov{d})\to (\ov{X}_*^{*-2},\ov{b},\ov{d})$.

\begin{proposition} The maps
\begin{align*}
& S_n\colon \HC_n(E,M)\to \HC_{n-2}(E,M),\\
& B_n\colon \HC_n(E,M)\to \HH_{n+1}(E,M),\\
& i_n\colon \HH_n(E,M)\to \HC_n(E,M),
\end{align*}
are induced by $-{}^e\ov{\varsigma}$, $i$ and $\pi$, respectively.\label{pr4.5}
\end{proposition}

\begin{proof} Left to the reader.
\end{proof}

Let $C$ be a $k$-algebra and let $M$ be a two-sided ideal of $C$. In \cite{G2} was proved that
if $M^{m+1}=0$, then
\[
S^{m(n+1)}\colon\HC_{n+2m(n+1)}(C,M)\to \HC_n(C,M)
\]
is the zero map. Actually, arguing as in the proof of Theorem~\ref{th4.6} it is easy to see
that if the previous formula holds when $m = 1$, then it is valid for all $m$ whenever
$M^{2^m} = 0$. Next, we give a theorem that improves this result.

\begin{theorem} If $M^{2^m} = 0$, then
\[
S^{m([n/2]+1)}\colon\HC_{n+2m([n/2]+1)}(C,M)\to \HC_n(C,M)\qquad\text{($n\ge 0$)}
\]
is the zero map, where $[n/2]$ denotes the integer part of $n/2$.\label{th4.6}
\end{theorem}

\begin{proof} We make the proof by induction on $m$. For $m=1$ the theorem follows from
Proposition~\ref{pr4.5}. Assume that $m>1$ and that the corollary is valid for $m-1$. Let $l =
[\frac{n}{2}]+1$. Consider the commutative diagram with exact rows
\[
\spreaddiagramcolumns{-1pc} \xymatrix{\HC_{n+2ml}(C,M^2) \rto\dto^{S^l} & \HC_{n+2ml}(C,M)
\rto^-{\pi} \dto^{S^l} &\HC_{n+2ml}\bigl(\frac{C}{M^2}, \frac{M}{M^2}\bigr)\dto^{S^l}\\
\HC_{n+2(m-1)l}(C,M^2) \rto^-{i}\dto^{S^{(m-1)l}} & \HC_{n+2(m-1)l}(C,M)\rto
\dto^{S^{(m-1)l}}& \HC_{n+2(m-1)l} \bigl(\frac{C}{M^2},\frac{M}{M^2}\bigr)\dto^{S^{(m-1)l}}\\
\HC_n(C,M^2) \rto^-{i} & \HC_n(C,M)\rto & \HC_n\bigl(\frac{C}{M^2},\frac{M}{M^2}\bigr)}
\]
where $i$ and $\pi$ are the canonical maps. Let $\bx\in \HC_{n+2ml}(C,M)$. By the case $m =
1$, we know that $S^l(\pi(\bx)) = 0$. So, there exists $\byy\in \HC_{n+2(m-1)l}(C,M^2)$ such
that $S^l(\bx) = i(\byy)$. Since $(M^2)^{2^{m-1}} = M^{2^m} = 0$ we can apply the inductive
hypothesis to conclude that $S^{ml}(\bx) = i(S^{(m-1)l}(S^l(\byy))) = 0$.
\end{proof}

Note that Theorem~\ref{th4.6} implies that if $M^2 = 0$, then the maps
\[
S_2\colon\HC_2(C,M)\to \HC_0(C,M) \quad\text{and}\quad S_3\colon\HC_3(C,M)\to \HC_1(C,M)
\]
are zero.

\subsection{The connection map for the Hochschild homology revised} Let
$$
\wh{\de}^3\colon A \ot \ov{A}^{\ot n}\to \ov{X}_{n-1}^1\quad\text{and}\quad\wt{\de}^3\colon
A\ot \ov{A}^{\ot n}\to \ov{X}_{n-1}^1
$$
be the maps defined by
\begin{align*}
&\wh{\de}^3(\ba)= \sum_{j=2}^{n-1}\frac{n\!+\!1\!-\!2j}{2(n\!+\!1)} F_{0j}(\ba) +
\sum_{i=0}^{n-1} \frac{i\!+\!1\!-\!n}{n\!+\!1} F_{in}(\ba) + \sum_{0<i<j<n}
\frac{i\!-\!j}{n\!+\!1} F_{ij}(\ba),\\
&\wt{\de}^3(\ba)=-\frac{1}{2} \sum_{0<i<j\le n} F_{ij}(\ba),
\end{align*}
where $F_{ij}(\ba)$ is the class of $t\xcirc F_i\xcirc F_j(\ba)$ in $\ov{X}_{n-1}^1$.

\begin{proposition} The connection map $\de_n\colon \HH_n(A)\to \HH_{n-1}(E,M)$ associated
with the short exact sequence~\eqref{eq3} is induced by the morphisms of complexes
\[
\wh{\de}\colon (A\ot\ov{A}^{\ot *},b)\to \Tot(\wt{X},\wt{\mathfrak{d}},\wt{\mathfrak{b}}
+{}^e\wt{\xi})\quad\text{and}\quad \wt{\de}\colon (A\ot\ov{A}^{\ot *},b)\to
\Tot(\wt{X},\wt{\mathfrak{d}}+{}^e\wt{\varsigma},\wt{\mathfrak{b}}),
\]
given by
\[
\wh{\de}(\ba) = (\de^2(\ba),\de^1(\ba),\wh{\de}^3(\ba),0,0,\dots)\quad\text{and} \quad
\wt{\de}(\ba) = (\de^2(\ba),\de^1(\ba),\wt{\de}^3(\ba),0,0,\dots),
\]
where we are writing
\begin{align*}
\Tot((\wt{X},\wt{\mathfrak{d}},\wt{\mathfrak{b}} +{}^e\wt{\xi}))_{n-1} &=
\Tot((\wt{X},\wt{\mathfrak{d}}+{}^e\wt{\varsigma},\wt{\mathfrak{b}}))_{n-1}\\
&=\wt{X}^{-1}_{n-2}\oplus \wt{X}^0_{n-1}\oplus \wt{X}^1_n\oplus\wt{X}^2_{n+1}\oplus
\wt{X}^3_{n+2}\oplus\cdots\\
& \simeq \ov{X}^0_{n-2}\oplus \ov{X}^0_{n-1}\oplus \ov{X}^1_{n-1}\oplus \ov{X}^1_n \oplus
\ov{X}^2_n\oplus \cdots
\end{align*}
(here we are identifying $\wt{X}^{-1}_{n-2} = \ov{X}^{-1}_{n-2}\oplus \ov{X}^0_{n-2}$ with
$\ov{X}^0_{n-2}$).\label{pr4.7}
\end{proposition}

\begin{proof} Let $\breve{\de}$ be the morphism introduced in Proposition~\ref{pr3.7}. By
Theorem~\ref{th4.3}, to check the assertions it suffices to show that $\Tot(\Psi)_{n-1}\xcirc
\wh{\de}_n = \Tot(P)_n \xcirc \breve{\de}_n$ and $\Tot(\Lambda)_{n-1}\xcirc\wh{\de}_n =
\wt{\de}_n$. In other words we must prove that
\begin{align}
&\Psi^{-1}(0,\de^2(\ba)) = (0,\de^2(\ba))\label{eq17},\\
&\Psi^0(\de^1(\ba),\wh{\de}_n^3(\ba)) = P(\de^1(\ba),\de^3(\ba)),\label{eq18}\\
&\Lambda^{-1}(0,\de^2(\ba)) = (0,\de^2(\ba)),\label{eq19}\\
&\Lambda^0(\de^1(\ba),\wh{\de}^3(\ba)) = (\de^1(\ba),\wt{\de}^3(\ba)),\label{eq20}
\end{align}
where $\ba\in A\ot \ov{A}^{\ot n}$ and $\de^1(\ba)$, $\de^2(\ba)$ and $\de^3(\ba)$ are as in
Subsection~\ref{TheConnection}. Equalities~\eqref{eq17} and~\eqref{eq19} are immediate. To
check the equality~\eqref{eq18} we first compute $\breve{P}(\de^1(\ba),\de^3(\ba))$. Since,
by~\eqref{eq15} we know that $(\de^1(\ba),-\si'\xcirc d\xcirc \de^1(\ba))\in P(\ddot{X})$ and
$P$ is a projection, we have:
\begin{align*}
P(\de^1(\ba),\de^3(\ba)) & = P(\de^1(\ba),-\si'\xcirc d\xcirc \de^1(\ba))+
P(0,(\de^3+ \si'\xcirc d\xcirc \de^1)(\ba))\\
&= \bigl(\de^1(\ba),-\si'\xcirc d\xcirc \de^1(\ba)\bigr)+ P\xcirc
\ddot{d}\sR\xcirc\bigl((\de^3+\si'\xcirc d\xcirc \de^1)(\ba),0\bigr)\\
& = \Bigl(\de^1(\ba),\frac{1}{2}\de^3(\ba)+\frac{1}{2}t\xcirc \de^3(\ba)
-\frac{1}{2}\si'\xcirc d\xcirc \de^1(\ba)+ \frac{1}{2}t\xcirc \si'\xcirc d\xcirc
\de^1(\ba)\Bigr).
\end{align*}
where the third equality follows from~\eqref{eq11}. Let
$$
L(\ba):= \frac{1}{2}\de^3(\ba)+\frac{1}{2}t\xcirc \de^3(\ba) -\frac{1}{2}\si'\xcirc d\xcirc
\de^1(\ba)+ \frac{1}{2}t\xcirc \si'\xcirc d\xcirc \de^1(\ba).
$$
By the definition of $\Psi^0$, to prove equality~\eqref{eq18}, we must show that
\begin{equation}
(-\ide + \ov{N}\xcirc \mathfrak{p})\xcirc \si'\xcirc d\xcirc \de^1(\ba)+
\ov{N}(\wh{\de}^3(\ba)) = L(\ba)\label{eq21}.
\end{equation}
A direct computation shows that
$$
t\xcirc F_i\xcirc F_j(\ba) = \begin{cases} F_{i+n-j}\xcirc t\xcirc F_j(\ba) &\text{if $0\le
i<j-1<n-1$,}\\ F_{i+1}\xcirc t\xcirc F_n(\ba) &\text{if $0\le i<n-2$ and $j=n$,}\end{cases}
$$
and
$$
t^2\xcirc F_i\xcirc F_j(\ba) = \begin{cases} - F_{j-i-1}\xcirc t\xcirc F_i(\ba) &\text{if
$0\le i<j-1<n-1$,}\\ - F_{n-i-2}\xcirc t\xcirc F_{i+1}(\ba) &\text{if $0\le i<n-2$ and $j=n$.}
\end{cases}
$$
Hence,
$$
\si'\xcirc d\xcirc \de^1(\ba) = \frac{1}{2}\sum_{i=1}^{n-2} \sum_{j=0}^n F_i\xcirc t\xcirc
F_j(\ba) = \frac{1}{2} \sum_{0\le i<j\le n} t\xcirc F_i\xcirc F_j(\ba) - \frac{1}{2}
\sum_{0\le i<j\le n} t^2\xcirc F_i\xcirc F_j(\ba).
$$
Thus, by the definition of $\de^3$,
$$
L(\ba) = \frac{1}{2} \sum_{j=2}^{n-1} F_0\xcirc F_j(\ba)+\frac{1}{2} \sum_{j=2}^{n-1} t\xcirc
F_0\xcirc F_j(\ba) - \sum_{0\le i<j\le n} t\xcirc F_i\xcirc F_j(\ba).
$$
and, by the definition of $\mathfrak{p}$,
\begin{align*}
\mathfrak{p}\xcirc\si'\xcirc d\xcirc\de^1(\ba) & =\sum_{0\le i<j< n}\frac{j-i}{2(n+1)}[t\xcirc
F_i\xcirc F_j(\ba)] + \sum_{i=0}^{n-1}
\frac{n-i-1}{2(n+1)}[t\xcirc F_i\xcirc F_n(\ba)]\\
&-\sum_{0\le i<j<n}\frac{n-j+i+1}{2(n+1)} [t^2\xcirc F_i\xcirc F_j(\ba)]
-\sum_{i=0}^{n-1}\frac{i+2}{2(n+1)}[t^2\xcirc F_i\xcirc F_n(\ba)]\\
& =\sum_{0\le i<j< n}\frac{2j-2i-n-1}{2(n+1)}[t\xcirc F_i\xcirc F_j(\ba)] + \sum_{i=0}^{n-1}
\frac{n-2i-3}{2(n+1)}[t\xcirc F_i\xcirc F_n(\ba)],
\end{align*}
where $[\bx]$ denotes the class of $\bx\in X_{n-1}^1$ in $\ov{X}_{n-1}^1$. Consequently,
\begin{align*}
& L(\ba)+(\ide - \ov{N}\xcirc \mathfrak{p})\xcirc \si'\xcirc d\xcirc \de^1(\ba) =
\sum_{j=2}^{n-1} \frac{n+1-2j}{2(n+1)} (t+t^2)\xcirc F_0\xcirc F_j(\ba)\\
&+ \sum_{i=0}^{n-1} \frac{i+1-n}{n+1} (t+t^2)\xcirc F_i\xcirc F_n(\ba)+ \sum_{0<i<j<n}
\frac{i-j}{n+1} (t+t^2)\xcirc F_i\xcirc F_j(\ba).
\end{align*}
The equality~\eqref{eq21} follows immediately from this fact. To prove equality~\eqref{eq20}
we must show that
$$
\wh{\de}_n^3(\ba) + \mathfrak{p}\xcirc \si' \xcirc d \xcirc \de^1(\ba) = -\frac{1}{2}
\sum_{0<i<j\le n} [t\xcirc F_i\xcirc F_j(\ba)],
$$
which can be checked by a direct computation.
\end{proof}

\appendix
\section{}
This appendix is devoted to prove Theorem~\ref{th3.2} and Propositions~\ref{pr3.7}
and~\ref{pr3.8}.

\smallskip

Let $E = A\ltimes_f M$ be a square zero extension. By definition, the Hochschild homology of
$E$ relative to $M$ is the homology of the complex
\[
(\breve{\mathfrak{X}},\breve{\mathfrak{b}}) = \xymatrix@1{{\ker\bigl((E\ot \ov{E}^{\ot *},b)}
\rto^-{\pi} & {(A\ot\ov{A}^{\ot *},b)\bigr)}},
\]
where $\pi$ is the canonical projection. Let $\hat{\mathfrak{X}}^w_v = (A\ot B^n_{w+1}) \oplus
(M\ot B^n_w)$. It is immediate that $(\breve{\mathfrak{X}},\breve{\mathfrak{b}})$ is the total
complex of the second quadrant double complex
\[
\spreaddiagramcolumns{0.75pc} \spreaddiagramrows{-0.25pc}
\xymatrix{\\\\\\
(\hat{\mathfrak{X}},\hat{\mathfrak{b}},\hat{\mathfrak{d}})=}\qquad\quad
\xymatrix{\vdots \dto^-{\hat{\mathfrak{b}}}&\vdots \dto^-{\hat{\mathfrak{b}}}&\vdots
\dto^-{\hat{\mathfrak{b}}}\\
\hat{\mathfrak{X}}^2_4&\hat{\mathfrak{X}}^1_4\dto^{\hat{\mathfrak{b}}}
\lto_{\hat{\mathfrak{d}}}& \hat{\mathfrak{X}}^0_4 \dto^{\hat{\mathfrak{b}}}
\lto_{\hat{\mathfrak{d}}}\\
& \hat{\mathfrak{X}}^1_3 \dto^{\hat{\mathfrak{b}}} & \hat{\mathfrak{X}}^0_3
\dto^{\hat{\mathfrak{b}}}\lto_{\hat{\mathfrak{d}}}\\
& \hat{\mathfrak{X}}^1_2 & \hat{\mathfrak{X}}^0_2 \dto^{\hat{\mathfrak{b}}}
\lto_{\hat{\mathfrak{d}}}\\
&& \hat{\mathfrak{X}}^0_1 \dto^{\hat{\mathfrak{b}}}\\
&& \hat{\mathfrak{X}}^0_0 }
\]
where $\hat{\mathfrak{d}}\colon \hat{\mathfrak{X}}^w_v \to \hat{\mathfrak{X}}^{w+1}_v$ is
defined by $\hat{\mathfrak{d}}(\bx_0^n) = \sum_{j=0}^{n-1}F_j(\bx_0^n)+ t\xcirc F_n(\bx_0^n)$
and $\hat{\mathfrak{b}}$ is given by the same formula as the Hochschild boundary map.

\smallskip

The $w$-th column $(\hat{\mathfrak{X}}^w,\hat{\mathfrak{b}})$ of the above complex is the
total complex of the double complex
\[
\spreaddiagramcolumns{1.2pc}
\xymatrix{\\\\
(\mathfrak{X}^w,\mathfrak{b},\alpha) =}\qquad\quad
\xymatrix{\vdots \dto^-{\mathfrak{b}_0} &\vdots \dto^-{\mathfrak{b}_1}\\
\mathfrak{X}^w_{0,2w+2}\dto^-{\mathfrak{b}_0}& \mathfrak{X}^w_{1,2w+2}
\dto^-{\mathfrak{b}_1}\lto_-{\alpha}\\
\mathfrak{X}^w_{0,2w+1}\dto^-{\mathfrak{b}_0}&\mathfrak{X}^w_{1,2w+1}
\dto^-{\mathfrak{b}_1} \lto_-{\alpha}\\
\mathfrak{X}^w_{0,2w} & \mathfrak{X}^w_{1,2w} \lto_-{\alpha}}
\]
where $\mathfrak{X}^w_{0v} = M\ot B_w^n$, $\mathfrak{X}^w_{1,v-1} = A\ot B_{w+1}^n$,
$\mathfrak{b}_0$ is given by the same formula as the Hochschild boundary map,
$\mathfrak{b}_1(\bx_0^n) = \mu_0^A(\bx_0^n)+\sum_{j=1}^{n-1}\mu_l(\bx_0^n)+ \mu_n^A(\bx_0^n)$
and $\alpha(\bx_0^n) = \mu_0^M(\bx_0^n) + \mu_n^M(\bx_0^n)$.

\begin{lemma} Let $\theta^w_1\colon (\mathfrak{X}^w_1,\mathfrak{b}_1)\to (X^w_1,-b)$ and
$\vartheta^w_1\colon (X^w_1,-b)\to (\mathfrak{X}^w_1,\mathfrak{b}_1)$ be the morphisms of
complexes given by
\[
\theta^w_1(\bx_0^{n+1}) = \mu_0^M(\bx_0^{n+1})\quad\text{and}\quad \vartheta^w_1(\bx_0^n) =
\sum_{l=0}^{n-i(\bx_0^n)} 1\ot\mathfrak{t}^l(\bx_0^n),
\]
where $\mathfrak{t}(\bx_0^n) = (-1)^n x_n\ot\bx_0^{n-1}$. Then, $\theta^w_1\xcirc
\vartheta^w_1 = \ide$ and $\vartheta^w_1\xcirc \theta^w_1$ is homotopic to $\ide$. A homotopy
is the family of maps $\ep^w\colon \mathfrak{X}^w_{1,v-1}\to \mathfrak{X}^w_{1v}$, defined by
\[
\ep^w(\bx_0^n)= -\sum_{l=0}^{n-i(\bx_0^n)} 1\ot \mathfrak{t}^l(\bx_0^n).
\]
\label{leA.1}
\end{lemma}

\begin{proof} It is immediate that $\theta^w_1$ is a morphism of chain complexes
and $\theta^w_1\xcirc \vartheta^w_1 = \ide$. We claim that $\vartheta^w_1$ is also a morphism
of chain complexes. Let $\bx\in X^w_{1v}$ be an elementary tensor and let $i=i(\bx)$. On one
hand, a direct computation shows that
$$
\mu_0^A(1\ot \bx) = \mu_{n+1}^A(1\ot \mathfrak{t}^{n-i}(\bx)) = 0\quad\text{and}\quad
\mu_0^A(1\ot \mathfrak{t}^l(\bx)) = -\mu_{n+1}^A(1\ot \mathfrak{t}^{l-1}(\bx)),
$$
for $1\le l\le n-i$, and so
$$
\mathfrak{b}_1(\vartheta^w_1(\bx)) = \sum_{j=1}^n \sum_{l=0}^{n-i}\mu_j(1\ot
\mathfrak{t}^l(\bx)) = - 1\ot \sum_{j=0}^{n-1}\sum_{l=0}^{n-i}\mu_j(\mathfrak{t}^l(\bx)).
$$
On the other hand, it is easy to see that
\begin{equation}
\mathfrak{t}^l(\mu_j(\bx)) =\begin{cases} \mu_{j+l}(\mathfrak{t}^l(\bx)) & \text{if $0\le l\le
n-i$ and $0\le j\le n-l-1$,}\\\mu_{j+l-n}(\mathfrak{t}^{l+1}(\bx)) & \text{if $0\le l\le
n-i-1$ and $n-l\le j\le n$,}\label{eq22}
\end{cases}
\end{equation}
and so
$$
\vartheta^w_1(b(\bx)) = \sum_{j=0}^{i-1}1\ot \mathfrak{t}^{n-i}(\mu_j(\bx))
+\sum_{l=0}^{n-i-1}\sum_{j=0}^n 1\ot \mathfrak{t}^l(\mu_j(\bx)) = 1\ot
\sum_{j=0}^{n-1}\sum_{l=0}^{n-i}\mu_j(\mathfrak{t}^l(\bx)),
$$
which proves the claim. We now check that $\ep^w$ is an homotopy from $\vartheta^w_1\xcirc
\theta^w_1$ to the identity map. Let $\bx\in \mathfrak{X}^w_{1,v-1}$ be an elementary tensor
and let $i = i(\bx)$. On one hand, since
$$
\mu_{n+1}^A(1\ot \mathfrak{t}^{n-i}(\bx)) = 0\quad\text{and}\quad \mu_0^A(1\ot
\mathfrak{t}^l(\bx)) = -\mu_{n+1}^A(1\ot \mathfrak{t}^{l-1}(\bx))\,\text{ for $1\le l\le
n-1$},
$$
we have
$$
\mathfrak{b}_1(\ep^w(\bx)) = -\bx - \sum_{j=1}^n\sum_{l=0}^{n-i}\mu_j(1\ot
\mathfrak{t}^l(\bx))= -\bx + \sum_{j=0}^{n-1}\sum_{l=0}^{n-i}1\ot \mu_j(\mathfrak{t}^l(\bx)).
$$
On the other hand, using the equality
$$
\sum_{l=0}^{n-i-1} 1\ot \mathfrak{t}^l\mu_n^A(\bx))=\sum_{l=0}^{n-i-1} 1\ot
\mathfrak{t}^l\mu_n(\bx))
$$
and the equations~\eqref{eq22} twice, we obtain
\begin{align*}
\ep^w((\mathfrak{b}_1(\bx)) & = -\sum_{l=0}^{n-i} 1\ot \mathfrak{t}^l(\mu_0^A(\bx)) -
\sum_{l=0}^{n-i-1} \sum_{j=1}^n 1\ot \mathfrak{t}^l(\mu_j(\bx)) - \sum_{j=1}^{i-1} 1\ot
\mathfrak{t}^{n-i}(\mu_j(\bx))\\
&= -\sum_{l=0}^{n-i} 1\ot \mathfrak{t}^l(\mu_0^A(\bx)) - \sum_{j=0}^{n-1} \sum_{l=0}^{n-i}
1\ot \mu_j(\mathfrak{t}^l(\bx)) + \sum_{l=0}^{n-i} 1\ot \mu_l(\mathfrak{t}^l(\bx)).
\end{align*}
Since $ \mathfrak{t}^l(\mu_0(\bx)) - \mu_l(\mathfrak{t}^l(\bx)) = \mathfrak{t}^l(\mu_0(\bx)) -
\mathfrak{t}^l(\mu^A_0(\bx))= \mathfrak{t}^l(\mu^M_0(\bx))$, we have
$$
\mathfrak{b}_1(\ep^w(\bx)) + \ep^w((\mathfrak{b}_1(\bx)) = -\bx + \sum_{l=0}^{n-i} 1\ot
\mathfrak{t}^l(\mu_0^M(\bx)) = -\bx + \vartheta^w_1(\theta^w_1(\bx)),
$$
which proves that $\ep^w$ is a homotopy from $\vartheta^w_1\xcirc \theta^w_1$ to $\ide$.
\end{proof}

\begin{lemma} For $w\ge 0$, let $\tau_0^1(\cX^w)$ be the double diagram with two columns
$\xymatrix{(X^w,b)&\lto_-{\ide-t} (X^w,-b)}$. The following assertions hold:

\begin{enumerate}

\smallskip

\item $\tau_0^1(\cX^w)$ is a double complex.

\smallskip

\item The map $\vartheta^w\colon \tau_0^1(\cX^w) \to (\mathfrak{X}^w,
\mathfrak{b},\alpha)$, where $\vartheta^w_{0}\colon X_v^w\to \mathfrak{X}_v^w$ is the identity
map and $\vartheta^w_{1}\colon X_{v-1}^w\to \mathfrak{X}_{v-1}^w$ is as in Lemma~\ref{leA.1},
is a morphism of double complexes.

\smallskip

\item The map $\hat{\theta}^w\colon (\hat{\mathfrak{X}}^w,\hat{\mathfrak{b}})\to
\Tot(\tau_0^1(\cX^w))$, defined by
\[
\qquad\qquad\hat{\theta}^w(\bx,\byy) = (\bx + t(\byy),\mu_0^M(\byy)),
\]
is a morphism of complexes.

\smallskip

\item Let $\hat{\vartheta}^w\colon \Tot(\tau_0^1(\cX^w)) \to (\hat{\mathfrak{X}}^w,
\hat{\mathfrak{b}})$ be the map induced by $\vartheta^w$. It is true that $\hat{\theta}^w
\xcirc \hat{\vartheta}^w = \ide$ and $\hat{\vartheta}^w\xcirc \hat{\theta}^w$ is homotopic to
the identity map. A homotopy is the family of maps
\[
\hat{\ep}^w\colon \mathfrak{X}^w_{0v}\oplus \mathfrak{X}^w_{1,v-1}\to \mathfrak{X}^w_{0,v+1}
\oplus \mathfrak{X}^w_{1v},
\]
defined by $\hat{\ep}^w(\bx,\byy)= (0,\ep^w(\byy))$, where $\ep^w$ is the homotopy introduced
in Lemma~\ref{leA.1}.
\end{enumerate}
\label{leA.2}
\end{lemma}

\begin{proof} By Lemma~\ref{leA.1} there is a special deformation retract
\[
\spreaddiagramcolumns{1.2pc} \xymatrix{{(X^w_*\!\oplus\! X^w_{*-1},b \!\oplus\!-b)}
\ar@<-1ex>[rr]_-{\ide \oplus \vartheta^w_1} && {(\mathfrak{X}^w_{0*}\!\oplus\!
\mathfrak{X}^w_{1,*-1},\mathfrak{b}_0\! \oplus\! \mathfrak{b}_1)} \ar@<-1ex>[ll]_-{\ide\oplus
\theta^w_1}},
\]
with homotopy $\hat{\ep}^w\colon \mathfrak{X}^w_{0*}\oplus \mathfrak{X}^w_{1,*-1}\to
\mathfrak{X}^w_{0,*+1}\oplus\mathfrak{X}^w_{1*}$, given by $\hat{\ep}^w(\bx,\byy) =
(0,\ep^w(\byy))$. Applying the perturbation lemma to this endowed with the perturbation
$\alpha$, we obtain a special deformation retract
\[
\spreaddiagramcolumns{1.2pc} \xymatrix{{\wh{\Tot}(\tau_0^1(\cX^w))}
\ar@<-1ex>[r]_-{\wt{\vartheta}^w} & {(\hat{\mathfrak{X}}^w,\hat{\mathfrak{b}})}
\ar@<-1ex>[l]_-{\wt{\theta}^w}},
\]
with homotopy $\wt{\ep}^w\colon \mathfrak{X}^w_{0*}\oplus \mathfrak{X}^w_{1,*-1} \to
\mathfrak{X}^w_{0,*+1}\oplus\mathfrak{X}^w_{1*}$. To finish the proof it remains to check that
\[\wh{\Tot}(\tau_0^1(\cX^w)) = \Tot(\tau_0^1(\cX^w)),\quad \wt{\vartheta}^w =
\hat{\vartheta}^w,\quad \wt{\theta}^w = \hat{\theta}^w\quad\text{and}\quad \wt{\ep}^w =
\hat{\ep}^w,
\]
which follow easily from the fact that $\ide - t = \alpha\xcirc\vartheta^w_1$ and
$\alpha\xcirc \ep^w = t$.
\end{proof}

The first item of the following lemma is part of item~(1) of Theorem~\ref{th3.2}.

\begin{lemma} The following assertions hold:

\smallskip

\begin{enumerate}

\item The diagram $(\hat{X},\hat{b},\hat{d})$, introduced above Theorem~\ref{th3.2}, is a
double complex.

\smallskip

\item The map $\hat{\vartheta}\colon (\hat{X},\hat{b},\hat{d})\to (\hat{\mathfrak{X}},
\hat{\mathfrak{b}},\hat{\mathfrak{d}})$, obtained by gluing the maps $\hat{\vartheta}^w$
introduced in Lemma~\ref{leA.2}, is a morphism of double complexes.

\smallskip

\item Recall from the proof of Lemma~\ref{th3.1} that $(\breve{X},\breve{b}) =
\Tot(\hat{X},\hat{b},\hat{d})$ and from the beginning of this section that
$(\breve{\mathfrak{X}},\breve{\mathfrak{b}}) = \Tot(\hat{\mathfrak{X}},
\hat{\mathfrak{b}},\hat{\mathfrak{d}})$. For $w\ge 0$ and $v\ge 2w$, let $ \hat{\zeta}^w_v
\colon \hat{\mathfrak{X}}^w_v\to \hat{X}^{w+1}_{v+1}$ be the maps defined by
\[
\hat{\zeta}^w(\bx_0^n,\byy_0^n) = \Biggl(0,F_0(\byy_0^n)+ \sum_{j=i(\byy_0^n)+1}^n
F_j(\byy_0^n)\Biggr).
\]
The map $\breve{\theta}\colon (\breve{\mathfrak{X}},\breve{\mathfrak{b}}) \to
(\breve{X},\breve{b})$, defined by $\breve{\theta}_n = \bigoplus_{w=0}^n \hat{\theta}^w_{n+w}
+ \bigoplus_{w=0}^{n-1} \hat{\zeta}^w_{n+w}$, where $\hat{\theta}^w_{n+w}$ is as in
Lemma~\ref{leA.2}, is a morphism of complexes.

\smallskip

\item Let $\breve{\vartheta}\colon (\breve{X},\breve{b})\to (\breve{\mathfrak{X}},
\breve{\mathfrak{b}})$ be the map induced by $\hat{\vartheta}$. It is true that
$\breve{\theta} \xcirc \breve{\vartheta} = \ide$ and $\breve{\vartheta}\xcirc \breve{\theta}$
is homotopic to the identity map. A homotopy is the family of maps $\breve{\ep}\colon
\breve{\mathfrak{X}}_n \to \breve{\mathfrak{X}}_{n+1}$, defined by $\breve{\ep}_{n+1} =
\bigoplus_{w=0}^n \hat{\ep}^w_{n+w+1}$, where $\hat{\ep}^w_{n+w+1}$ is as Lemma~\ref{leA.2}.
\end{enumerate}
\label{leA.3}
\end{lemma}

\begin{proof} By Lemma~\ref{leA.2} we have the following special deformation retract:
\[
\xymatrix{{}\save[]+<-34pt,-4pt>\Drop{\displaystyle{\bigoplus_{w\ge 0}}
(\hat{X}^w_{*+w},\hat{b})}\restore \ar@<-1ex>[rr]_-{\hat{\vartheta}} &&
{{}\save[]+<34pt,-4pt>\Drop{\displaystyle{\bigoplus_{w\ge 0}} (\hat{\mathfrak{X}}^w_{*+w},
\hat{\mathfrak{b}})}} \restore \ar@<-1ex>[ll]_-{\hat{\theta}}},\qquad \hat{\ep},
\]
where $\hat{\vartheta} = \bigoplus_{w\ge 0} \hat{\vartheta}^w_{*+w}$, $\hat{\theta} =
\bigoplus_{w\ge 0} \hat{\theta}^w_{*+w}$ and $\hat{\ep} = \bigoplus_{w\ge 0}
\hat{\ep}^w_{*+w}$. Consider the perturbation $\hat{\mathfrak{d}}$. Applying the perturbation
lemma to this datum, we obtain a special deformation retract
\[
\xymatrix{{\ov{\Tot}(\hat{X},\hat{b},\hat{d})} \ar@<1ex>[r]^-{\ov{\vartheta}} &
{(\breve{\mathfrak{X}},\breve{\mathfrak{b}})} \ar@<1ex>[l]^-{\ov{\theta}}},\quad
\ov{\ep}_{*+1}\colon \breve{\mathfrak{X}}_*\to \breve{\mathfrak{X}}_{*+1}.
\]
To finish the proof it remains to check that
\[
\ov{\Tot}(\hat{X},\hat{b},\hat{d}) = (\breve{X},\breve{b}),\quad \ov{\vartheta} =
\breve{\vartheta},\quad \ov{\theta} = \breve{\theta}\quad\text{and}\quad \ov{\ep} =
\breve{\ep},
\]
for which it suffices to check that
$$
\hat{\theta}\xcirc \hat{\mathfrak{d}}\xcirc \hat{\vartheta} = \hat{d},\quad \hat{\ep}\xcirc
\hat{\mathfrak{d}}\xcirc \hat{\vartheta}= 0,\quad \hat{\theta}\xcirc \hat{\mathfrak{d}}\xcirc
\hat{\ep}= \hat{\zeta}\quad\text{and}\quad \hat{\ep}\xcirc \hat{\mathfrak{d}}\xcirc \hat{\ep}
= 0,
$$
where $\hat{\zeta}_n = \bigoplus_{w=0}^{n-1} \hat{\zeta}^w_{n+w}$, which follows by a direct
computation.
\end{proof}

Let $(\breve{\mathfrak{X}},\breve{\mathfrak{b}},\breve{\mathfrak{B}}) =
\xymatrix@1{{\ker\bigl((E\ot \ov{E}^{\ot *},b,B)} \rto^-{\pi} & {(A\ot\ov{A}^{\ot
*},b,B)\bigr)}}$.

\begin{lemma} Let $\hat{B}$ be as in Theorem~\ref{th3.2}. The following assertions hold:

\begin{enumerate}

\smallskip

\item $(\hat{X},\hat{b},\hat{d},\hat{B})$ is a double mixed complex.

\smallskip

\item Let $(\breve{X},\breve{b},\breve{B})$ be the mixed complex associated with
$(\hat{X},\hat{b},\hat{d},\hat{B})$. The maps
\[\qquad\quad
\breve{\vartheta}\colon (\breve{X},\breve{b},\breve{B}) \to (\breve{\mathfrak{X}},
\breve{\mathfrak{b}},\breve{\mathfrak{B}}) \quad \text{and}\quad \breve{\theta}\colon
(\breve{\mathfrak{X}},\breve{\mathfrak{b}},\breve{\mathfrak{B}}) \to (\breve{X},\breve{b},
\breve{B}),
\]
introduced in Lemma~\ref{leA.3}, are morphisms of mixed complexes.
\end{enumerate}
\label{leA.4}
\end{lemma}

\begin{proof} 1)\enspace From the fact that $\hat{b}\xcirc \hat{b} = 0$ it follows easily that
$t\xcirc b = b\xcirc t$. Thus we obtain that $b\xcirc N = N\xcirc b$, which implies that
$\hat{b}\xcirc \hat{B} + \hat{B}\xcirc \hat{b} = 0$. To prove that $\hat{d}\xcirc \hat{B} +
\hat{B}\xcirc \hat{d} = 0$ we must check that $d'\xcirc N = -N\xcirc d$. Let $\bx_0^n\in
X_v^w$ be an elementary tensor. Let $0=i_0<i_1<\cdots<i_w\le n$ be the indices such that
$x_{i_j}\in M$ and let $i_{w+1} = n+1$. Fix $l\in \{0,\dots,n-1\}\setminus(\{i_0,\dots,i_w\}
\cup \{i_1-1,\dots,i_w-1\})$ and let $r$ such that $i_r<l<i_{r+1}$. A direct computation shows
that
$$
t^j(F_l(\bx_0^n)) = \begin{cases} F_{l+n+1-i_{w+1-j}}\xcirc t^j(\bx_0^n)&\text{if $0\le j\le
w-r$,}\\ t\xcirc F_{l+n+1-i_{r+1}}\xcirc t^{w-r}(\bx_0^n) &\text{if $j = w-r+1$,}\\
F_{l-i_{w+2-j}}\xcirc t^{j-1}(\bx_0^n)&\text{if $j > w-r+1$.}
\end{cases}
$$
Hence
$$
d'\xcirc N(\bx_0^n) = -\sum_{l=1}^{n-1}\sum_{j=0}^w F_l\xcirc t^j(\bx_0^n) -
\sum_{l=i(\bx_0^n)+1}^{n-1}\sum_{j=0}^w t\xcirc F_l\xcirc t^j(\bx_0^n) = - \sum_{j=0}^{w+1}
\sum_{l=1}^{n-1} t^j\xcirc F_l(\bx_0^n),
$$
as we want.

\smallskip

\noi 2)\enspace From Lemma~\ref{leA.3} we get a special deformation retract between the total
complexes of the double complexes $\BC(\breve{X},\breve{b},0)$ and $\BC(\breve{\mathfrak{X}},
\breve{\mathfrak{b}},0)$. Consider the perturbation $\breve{\mathfrak{B}}$. The result it
follows by applying the perturbation lemma to this setting, and using that $\breve{B} =
\breve{\theta}\xcirc \breve{\mathfrak{B}}\xcirc \breve{\vartheta}$, $\breve{\mathfrak{B}}
\xcirc \breve{\ep} = 0$ and $\breve{\ep}\xcirc \breve{\mathfrak{B}} = 0$.
\end{proof}

\noi{\bf Proof of Theorem~\ref{th3.2}.}\enspace It follows immediately from
Lemma~\ref{leA.4}.\qed

\bigskip

\noi{\bf Proof of Proposition~\ref{pr3.7}.}\enspace It is immediate that the sequence
\[
\xymatrix{0\rto &{(\breve{\mathfrak{X}},\breve{\mathfrak{b}})} \rto^-{i} &{(E\ot\ov{E}^{\ot
*},b)}\rto^-{\pi} &{(A\ot\ov{A}^{\ot *},b)}\rto & 0.}
\]
splits in each level via the maps $s_n\colon A\ot \ov{A}^{\ot n} \to E\ot\ov{E}^{\ot n}$ and
$r_n\colon E\ot \ov{E}^{\ot n} \to \breve{\mathfrak{X}}_n$, given by $s_n(\ba_0^n) =
(a_0,0)\ot\cdots\ot (a_n,0)$ and $r_n = \ide - s_n\xcirc \pi_n$. From this it follows that the
connection map $\de_n\colon \HH_n(A)\to \HH_{n-1}(E,M)$ is induced by the map $\de\colon (A\ot
\ov{A}^{\ot *},b)\to (\breve{\mathfrak{X}},\breve{\mathfrak{b}})[1]$, given by $\de_n =
r_{n-1}\xcirc b_n\xcirc s_n$. To finish the proof it suffices to check that $\breve{\de} =
\breve{\theta}\xcirc \de$, where $\breve{\theta}\colon
(\breve{\mathfrak{X}},\breve{\mathfrak{b}}) \to (\breve{X},\breve{b})$ is as in
Lemma~\ref{leA.3}.\qed

\bigskip

\noi{\bf Proof of Proposition~\ref{pr3.8}.}\enspace From Proposition~\ref{pr3.7} and
Lemma~\ref{leA.4} it follows that $\de_n$ is induced by the morphism of complexes
\[
\acute{\de}\colon \Tot(\BC(A\ot \ov{A}^{\ot *},b,B))\to \Tot(\BC(\breve{X},\breve{b},
\breve{B}))[1],
\]
given by $\acute{\de}_n(\ba_0^n,\bb_0^{n-2},\bc_0^{n-4},\dots) =
(\breve{\de}_n(\ba_0^n),\breve{\de}_{n-2}(\bb_0^{n-2}),\breve{\de}_{n-4}(\bc_0^{n-4}),\dots)$,
where $\breve{\de}_n$, $\breve{\de}_{n-2}$, $\breve{\de}_{n-4}$, etcetera, are as in
Proposition~\ref{pr3.7}. To finish the proof it suffices to compose this map with the
canonical projection from $\Tot(\BC(\breve{X},\breve{b},\breve{B}))[1]$ to
$\Tot(\ov{X}_0,\ov{b}_0,\ov{d}_0)[1]$.\qed

\appendix
\setcounter{section}{1}
\section{}
Recall from the discussion above Theorem~\ref{th4.3}, that for each $t$-invariant element
$\bx\in X_v^w$,
\[
\Upsilon(\bx) = (-\ide + \ov{N}\xcirc \mathfrak{p})\xcirc \si'\xcirc d(\bx)\in X_v^{w+1}.
\]
It is easy to check that $\Upsilon(\bx)$ is univocally determined by the following properties:
\[
\Upsilon(\bx)\in \ker(\mathfrak{p})\qquad\text{and}\qquad (x,\Upsilon(\bx))\in P(\ddot{X}).
\]
Let ${}^e\! P(\ddot{X}_v^w) = \{(\bx,\Upsilon(\bx)) \in \ddot{X}_v^w: \bx \text{ is
$t$-invariant} \}$. Clearly, $P(\ddot{X}) = {}^e\! P(\ddot{X}) \oplus {}^t\!P(\ddot{X})$. We
assert that $\ddot{d}({}^e\! P(\ddot{X}_v^w)) \sub {}^e\! P(\ddot{X}_v^{w+1})$. In order to
prove this we will need the following result

\begin{lemma} It is true that $\ov{d}\xcirc \mathfrak{p}=-\mathfrak{p}\xcirc d'$.\label{leB.1}
\end{lemma}

\begin{proof} Let $\byy_0^n\in X_v^w$ and let $i_w>0$ be the last index
such that $y_{i_w}\in M$. We have: \allowdisplaybreaks
\begin{align*}
\mathfrak{p} \xcirc d'(\byy_0^n) & = -\sum_{j=1}^{n-1} \mathfrak{p}\xcirc F_j(\byy_0^n) -
\sum_{j=i_w+1}^{n-1} \mathfrak{p}\xcirc t \xcirc F_j(\byy_0^n)\\
& = - \sum_{j=1}^{i_w-2} \frac{n-i_w+2}{v+2} [F_j(\byy_0^n)] - \sum_{j=i_w+1}^{n-1}
\frac{n-j+1}{v+2}[F_j(\byy_0^n)]\\
& - \sum_{j=i_w+1}^{n-1} \frac{j-i_w+1}{v+2}[t\xcirc F_j(\byy_0^n)]\\
& = - \sum_{j=1}^{n-1} \frac{n-i_w+2}{v+2} [F_j(\byy_0^n)]\\
& = -\ov{d}\xcirc \mathfrak{p}(\byy_0^n),
\end{align*}
where $[\bx]$ denotes the class of $\bx\in X_v^{w+1}$ in $\ov{X}_v^{w+1}$.
\end{proof}

\begin{proposition} Let $\bx\in X_v^w$ be a $t$-invariant element. Then
\[
\ddot{d}(\bx,\Upsilon(\bx)) = -\frac{w+1}{w+2}(d'(\bx),\Upsilon\xcirc d'(\bx)).
\]
\label{prB.2}
\end{proposition}

\begin{proof} Since, by Lemma~\ref{le3.3},
$$
(t-\ide)(\Upsilon(\bx)) = d(\bx) - \frac{1}{w+2}N\xcirc d(\bx)= d(\bx) + \frac{1}{w+2}
d'\xcirc N(\bx) = d(\bx) + \frac{w+1}{w+2} d'(\bx),
$$
we have:
\[
\ddot{d}(\bx,\Upsilon(\bx)) = (d(\bx)+(\ide-t)(\Upsilon(\bx)),d'(\Upsilon(\bx)) =
\left(-\frac{w+1}{w+2} d'(\bx),d'(\Upsilon(\bx))\right).
\]
In order to finish the proof it suffices to check that
\[
\left(-\frac{w+1}{w+2} d'(\bx),d'(\Upsilon(\bx))\right)\in P(\ddot{X})\quad\text{and}\quad
\mathfrak{p}(d'(\Upsilon(\bx)) = 0.
\]
The first fact follows immediately from the fact that $P(\ddot{X})$ is a subcomplex of
$(\ddot{X},\ddot{b},\ddot{d})$ and the second one follows easily from Lemma~\ref{leB.1}.
\end{proof}

For each $v$ and $w$, let ${}^e b\colon {}^e\! P(\ddot{X}_v^w)\to {}^e\! P(\ddot{X}_{v-1}^w)$
and ${}^e\xi\colon {}^e\! P(\ddot{X}_v^w)\to {}^t\! P(\ddot{X}_{v-1}^w)$ be the maps defined
by $\ddot{b}(\bx) = {}^e b(\bx) + {}^e\xi(\bx)$. We now want to compute these maps. To carry
out this task we will need Proposition~\ref{prB.3} below. Let ${}^t\! X_v^w$ be the set of
$t$-invariant elements of $X_v^w$ and let $\xi\colon {}^t\! X_v^w\to {}^t\! X_{v-1}^{w+1}$ be
the map defined by
\[
\xi(\bx) = -\ov{N}\xcirc \mathfrak{p}\xcirc \si'\xcirc d\xcirc b(\bx) - \ov{N}\xcirc
\ov{b}\xcirc \mathfrak{p}\xcirc \si'\xcirc d(\bx).
\]

\begin{proposition} Assume that $\bx\in X_v^w$ is a $t$-invariant element. Then,
\[
{}^e b(\bx,\Upsilon(\bx)) = (b(\bx),\Upsilon(b(\bx)))\quad \text{and}\quad {}^e
\xi(\bx,\Upsilon(\bx)) = (0,\xi(\bx)).
\]
\label{prB.3}
\end{proposition}

\begin{proof} First note that $\ddot{b}(\bx,\Upsilon(\bx)) \in P(\ddot{X})$, since
$P(\ddot{X})$ is a subcomplex of $(\ddot{X},\ddot{b},\ddot{d})$. So, from the fact that
\[
\ddot{b}(\bx,\Upsilon(\bx)) = (b(\bx),-b(\Upsilon(\bx)) =(b(\bx),\Upsilon(b(\bx))) +(0,
-b(\Upsilon(\bx))-\Upsilon(b(\bx))),
\]
it follows that
\[
{}^e b(\bx,\Upsilon(\bx)) = (b(\bx),\Upsilon(b(\bx)))\quad\text{and}\quad
{}^e\xi(\bx,\Upsilon(\bx)) = (0,-\Upsilon(b(\bx)) - b(\Upsilon(\bx))).
\]
To finishes the proof we must compute the last map. But, since by Proposition~\ref{pr3.5},
\[
b\xcirc \si'\xcirc d(\bx) = \si'\xcirc b\xcirc d(\bx) = - \si'\xcirc d\xcirc b(\bx),
\]
we have:
\begin{align*}
-\Upsilon(b(\bx))-b(\Upsilon(\bx)) & = (\ide - \ov{N}\xcirc \mathfrak{p})\xcirc \si'\xcirc
d\xcirc b(\bx) + b\xcirc (\ide - \ov{N}\xcirc \mathfrak{p})\xcirc \si'\xcirc d(\bx)\\
& = - \ov{N}\xcirc \mathfrak{p}\xcirc \si'\xcirc d\xcirc b(\bx) - b\xcirc \ov{N}\xcirc
\mathfrak{p}\xcirc \si'\xcirc d(\bx)\\
& =-\ov{N}\xcirc \mathfrak{p}\xcirc \si'\xcirc d\xcirc b(\bx) - \ov{N}\xcirc \ov{b}\xcirc
\mathfrak{p}\xcirc \si'\xcirc d(\bx),
\end{align*}
as we want.
\end{proof}

\noi{\bf Proof of Theorem~\ref{th4.3}.}\enspace It is immediate that $\Lambda$ is a morphism
of double mixed complexes. Using Propositions~\ref{prB.2} and~\ref{prB.3} it is easy to check
that $\Psi$ is a morphism of double mixed complexes. To finish the proof it suffices to note
that the maps $\Lambda$ and $\Psi$ are bijective.\qed

\smallskip

We now are going to prove Proposition~\ref{pr4.4}. To carry out this task it is convenient to
first calculate $\mathfrak{p}\xcirc \si'\xcirc d\xcirc \ov{N}$.

\begin{lemma} Let $\bx_0^n\in X_v^w$ be an elementary tensor, let
$0=i_0<\dots<i_w\le n$ be the indices such that $x_{i_j}\in M$ and let $i_{w+1} = n+1$. Given
$0\le \al \le n$, we let $j(\al)$ denote the number defined by $i_{j(\al)}\le \al\le
i_{j(\al)+1}$. We have:
\[
\mathfrak{p}\xcirc \si'\xcirc d\xcirc \ov{N}([\bx_0^n]) = \sum_{\al=1}^{n-1} (-1)^{\al}
D_{\al} [F_{\al }(\bx_0^n)],
\]
where $[\bx_0^n]$ and $[F_{\al }(\bx_0^n)]$ denote the class of $\bx_0^n$ in $\ov{X}_v^w$ and
$F_{\al }(\bx_0^n)$ in $\ov{X}_v^{w+1}$, respectively, and
\[
D_{\al} = \frac{w+1}{2}-\frac{(w+1)(w+2\al+2) + 2(n+1)(w-j(\al))-2\sum_{u=1}^w i_u}
{2(w+2)(v+2)}.
\]
\label{leB.4}
\end{lemma}

\begin{proof} Let $d'' = -d' - d$. Then,
\[
d\xcirc \ov{N}([\bx_0^n]) = - d'\xcirc \ov{N}([\bx_0^n]) - d''\xcirc \ov{N}([\bx_0^n]) =
\ov{N}\xcirc \ov{d}([\bx_0^n])-d''\xcirc\ov{N}([\bx_0^n]).
\]
On one hand, it is immediate that
\[
\mathfrak{p}\xcirc \si' \xcirc \ov{N}\xcirc \ov{d}([\bx_0^n]) = \left(\sum_{j=0}^w
\frac{w+1-j}{w+2}\right) \mathfrak{p} \xcirc \ov{N} \xcirc \ov{d}([\bx_0^n])= \frac{w+1}{2}
\ov{d}([\bx_0^n]).
\]
On the other hand,
$$
\mathfrak{p}\xcirc \si'\xcirc d''\xcirc \ov{N}([\bx_0^n]) = \sum_{\al=1}^{n-1}
\mathfrak{p}\xcirc \si'\xcirc t\xcirc F_{\al}(\bx_0^n)= \sum_{\al=1}^{n-1} C_{\al}
[F_{\al}(\bx_0^n)],
$$
where
\begin{align*}
C_{\al} & = \frac{(w+1)(\al+1-i_{j(\al)})}{(w+2)(v+2)}+ \sum_{u=1}^{j(\al)}
\frac{(u+w-j(\al))(i_u-i_{u-1}+1)}{(w+2)(v+2)}\\
& + \frac{(w-j(\al))(n-i_w+2)}{(w+2)(v+2)} + \sum_{u=j(\al)+2}^{w}
\frac{(u-j(\al)-1)(i_u-i_{u-1}+1)} {(w+2)(v+2)}\\
& = \frac{(w+1)(w+2\al+2) + 2(n+1)(w-j(\al))-2\sum_{u=1}^w i_u}{2(w+2)(v+2)}.
\end{align*}
The result follows immediately from these facts.
\end{proof}

\noi{\bf Proof of Proposition~\ref{pr4.4}.}\enspace By Lemma~\ref{leB.4},
$$
\ov{d}\xcirc \mathfrak{p}\xcirc\si'\xcirc d\xcirc \ov{N}([\bx_0^n])\! =\! \sum_{\al<\be}
L_{\al\be} [F_{\al}\xcirc F_{\be}(\bx_0^n)],
$$
where
$$
L_{\al\be} = D_{\be}-D_{\al} = \frac{(w+1)(\al-\be)+(n+1)(j(\be)-j(\al))}{(w+2)(v+2)},
$$
and
$$
\mathfrak{p}\xcirc\si'\xcirc d\xcirc \ov{N}\xcirc \ov{d} ([\bx_0^n]) = \sum_{\al<\be}
L'_{\al\be} [F_{\al}\xcirc F_{\be}(\bx_0^n)],
$$
where
\begin{align*}
L'_{\al\be} & = \biggl(\frac{(w+2)(w+2\be+1) + 2n(w-j(\be))-2(\sum i_u-w+j(\al)+\al)}
{2(w+3)(v+2)}\\
& - \frac{(w+2)(w+2\al+3) + 2n(w+1-j(\al)) -2(\sum i_u-w+j(\be)+2\be)}{2(w+3)(v+2)}\biggr)\\
& = \frac{(w+3)(\be-\al) + (n-1)(j(\al)-j(\be))-(v+2)}{(w+3)(v+2)}.
\end{align*}
So,
$$
{}^e\wt{\varsigma}([\bx_0^n]) = \frac{1}{w+1}\ov{d}\xcirc\mathfrak{p}\xcirc \si'\xcirc d\xcirc
\ov{N}([\bx_0^n]) + \frac{1}{w+2}\mathfrak{p}\xcirc \si'\xcirc d\xcirc\ov{N} \xcirc
\ov{d}([\bx_0^n])  = \sum_{\al<\be}\la_{\al\be}^{(w)} [F_{\al}\xcirc F_{\be}(\bx_0^n)],
$$
where
\[
\la_{\al\be}^{(w)} =\frac{2(j(\be)-j(\al))}{(w+1)(w+2)(w+3)}- \frac{1}{(w+2)(w+3)},
\]
as desired.\qed


\begin{thebibliography}{G-S}

\bibitem[B]{B} D. Burghelea
\newblock {\em Cyclic homology and algebraic $K$-theory of spaces I},
\newblock {Boulder Colorado 1983, Contemp. Math., vol~55}
\newblock {(1986) 89--115}.

\bibitem[C-Q]{C-Q} J. Cuntz and D. Quillen
\newblock {\em Operators on noncommutative differential forms and cyclic homology},
\newblock {Geometry, Topology and Phisics; for Raoul Bott}
\newblock {International Press, Cambridge MA (1995) 77--111}.

\bibitem[Co]{Co} G. Corti\~nas
\newblock {\em On the cyclic homology of commutative algebras over arbitrary ground rings},
\newblock {Communications in Algebra, vol. 27:3}
\newblock {(1999) 1403--1412}.

\bibitem[C]{C} M. Crainic
\newblock {\em On the perturbation lemma, and deformations},
\newblock {arXiv:Math. AT/0403266}
\newblock {(2004)}.

\bibitem[G-S]{G-S} M. Gerstenhaber and S. D. Schack
\newblock {\em Relative Hochschild cohomology, rigid algebras and the Bockstein},
\newblock {Journal Of Pure and Applied Algebra, vol~43}
\newblock {(1986) 53--74}.

\bibitem[G1]{G1} T. G. Goodwillie
\newblock {\em Relative algebraic $K$-theory and cyclic homology},
\newblock {Ann. of Math., vol~124}
\newblock {(1986) 347--402}.

\bibitem[G2]{G2} T. G. Goodwillie
\newblock {\em  Cyclic homology, derivations and the free loop space},
\newblock {Topology, vol~24}
\newblock {(1985) 187--215}.

\bibitem[K1]{K1} L. Kadison
\newblock {\em  Cyclic homology of extension algebras with application to matrix algebras,
algebraic $K$-theory, and Nest algebras of operators},
\newblock {Ph. D. thesis, U. of Cal. (Berkeley)}
\newblock {(1984)}.

\bibitem[K2]{K2} L. Kadison
\newblock {\em A relative cyclic cohomology theory useful for computations},
\newblock {Cr. Acad. Sci. Paris, vol~308}
\newblock {(1989) 569--573}.

\bibitem[Ka1]{Ka1} K. Kassel
\newblock {\em Cyclic homology, comodules and mixed complexes},
\newblock {Journal of Algebra, vol~107}
\newblock {(1987) 195--216}.

\bibitem[Q]{Q} D. Quillen
\newblock {\em Cyclic cohomology and algebra extensions},
\newblock {K-theory, vol~3}
\newblock {(1989) 205--246}.

\end{thebibliography}
\end{document}